
\input amstex
\documentstyle{amsppt}
\magnification=\magstep1
\hoffset .2 true in \voffset .2 true in
\hsize=6.5 true in \vsize=8.6 true in

\def\Z{{\Bbb Z}}
\def\C{{\Bbb C}}
\def\A{{\Cal A}}
\def\E{{\Cal E}}
\def\dd{{\Cal D}}

\def\l{{\ell}}
\def\i{{\text{i}}}
\def\D{{\Delta}}
\def\DS{\displaystyle}
\def\id{\operatorname{id}}
\def\GL{\operatorname{GL}}
\def\PB{\operatorname{PB}}
\def\Aut{\operatorname{Aut}}
\def\IA{\operatorname{IA}}
\def\End{\operatorname{End}}
\def\Hom{\operatorname{Hom}}
\def\Homeo{\operatorname{Homeo}}
\def\rank{\operatorname{rank}}
\def\im{\operatorname{im}}
\def\ev{\operatorname{ev}}

\def\ab{\operatorname{ab}}

\topmatter

\title Homology of Iterated Semidirect Products of Free Groups
\endtitle

\rightheadtext{Iterated Semidirect Products of Free Groups}
\leftheadtext{Daniel C.~Cohen and Alexander I.~Suciu}

\author Daniel C.~Cohen$^1$ and Alexander I.~Suciu$^2$
\endauthor

\affil
\ \\
\line{$^1$\ Department of Mathematics,
University of California,
Davis, CA 95616 USA\hfil}\\
\line{$^{\ }$\ {\eightrm cohen\@math.ucdavis.edu}\hfil}\\
\\
\line{$^2$\ Department of Mathematics,
Northeastern University,
Boston, MA 02115 USA\hfil}\\
\line{$^{\ }$\ {\eightrm alexsuciu\@neu.edu}\hfil}\\
\endaffil

\thanks{The second author was partially supported by N.S.F.~grant
DMS--9103556, and an internal grant from the Mathematics Department
at Northeastern University.} \endthanks

\endsubjclass

\abstract
Let $G$ be a group which admits the structure of an iterated
semidirect product of finitely generated free groups.
We construct a finite, free resolution of the integers
over the group ring of $G$.  This resolution is used to
define representations of groups which act compatibly on $G$,
generalizing classical constructions of Magnus, Burau, and Gassner.
Our construction also yields algorithms for computing the homology
of the Milnor fiber of a fiber-type hyperplane arrangement,
and more generally, the homology of the complement of such
an arrangement with coefficients in an arbitrary local system.
\endabstract

\endtopmatter

\document

\head Introduction
\endhead

Let $G=F_{d_\l}\rtimes\cdots\rtimes F_{d_2}\rtimes F_{d_1}$ be a
group which admits the structure of an iterated semidirect product
of finitely generated free groups.  For any such group,
we construct an explicit finite, free resolution
$C_\bullet(G)\to\Z$ over the group ring of $G$
(Theorem~2.11). Topologically, this resolution may
be viewed as the equivariant chain complex of the universal cover of an
Eilenberg-MacLane space of type $K(G,1)$.
The boundary maps of the chain complex $C_\bullet(G)$ are computed
recursively by means of Fox derivatives from the various actions of
$F_{d_p}$ on $F_{d_q}$, $p<q$, dictated by the semidirect product
structure of $G$.  Independent of these actions, each term, $C_k(G)$, of
$C_\bullet(G)$ is a free $\Z G$-module of rank $\sum d_{p_1}d_{p_2}\cdots
d_{p_k}$ (the sum being over all $1 \le p_1 < \cdots < p_k \le \l$).

Perhaps the most famous groups of this type are Artin's pure braid
groups.  The pure braid group on $\l$ strings
may be realized as
$P_\l= F_{\l-1}\rtimes\cdots\rtimes F_2\rtimes F_1$, \cite{4}.
A natural generalization also belongs to this
class of groups.  The fundamental group of the complement of any
(affine) fiber-type hyperplane arrangement admits the structure
of an iterated semidirect product of free groups,
\cite{18}, \cite{25}, \cite{43}.  Examples include Brieskorn's
generalized pure braid groups, $\PB(W)$, where $W$ is a Coxeter group
of type $\text{A}_\l$, $\text{B}_\l$, $\text{G}_2$, or $\text{I}_2(p)$,
see~\cite{8}.  Groups of the form $P_{n,\l}=\ker (P_n\to P_\l)$
also admit this structure.  These latter groups arise in
the studies of representations of braid groups, generalized
hypergeometric functions, and the Knizhnik-Zamolodchikov equations,
see e.g.~\cite{29}, \cite{2}, and \cite{47}.

Much is known about many groups of this type.  Arnol'd \cite{3}
and Cohen \cite{13} computed the cohomology of the pure braid
group $P_\l = \PB (\text{A}_{\l-1})$, and showed that the Poincar\'e
polynomial factors into linear terms.  The lower central series of
$P_\l$ was found by Kohno~\cite{27}.  These two results combine to
yield the ``LCS formula'' relating the Betti numbers, $b_j$, the
exponents, $d_q=q$, and the ranks of the lower central series
quotients, $\phi_k$, of $P_\l$:
$$\sum_{j=0}^{\l-1} b_j(-t)^j = \prod_{q=1}^{\l-1}(1-d_qt) =
\prod_{k\ge 1}(1-t^k)^{\phi_k}.$$

These results were subsequently generalized to the group of an
arbitrary fiber-type arrangement by Falk and Randell \cite{18},
\cite{19} and Jambu \cite{25}. See \cite{30} for analogous results
on the pure braid group of a Riemann surface, and see \cite{8},
\cite{41}, \cite{42}, and \cite{48} for further results on the
cohomology and factorization of the Poincar\'e polynomial of an
arrangement.   We obtain a further generalization here.
The LCS formula holds for any iterated
semidirect product of free groups $G$ for which the split extensions
arising in the semidirect product structure give rise to
IA-automorphisms of the free groups comprising $G$ (Theorem~3.7).

The aforementioned results on the homology of an iterated
semidirect product of free groups, $G$, apply only in the constant
coefficient case (that is, homology with coefficients in a
trivial $G$-module).  A desire to compute the homology of $G$
with coefficients in an arbitrary $G$-module led to the construction
of this paper.  Let $\nu:G\to\Aut(V)$ be a representation of
$G$, and denote by $V={V}_{\nu}$ the corresponding $\Z
G$-module. Then the homology of $G$ with coefficients in $V$ is
equal to the homology of the chain complex $C_\bullet(G) \otimes_G
V$ (see \cite{9}).  In this manner, we obtain an algorithm for
computing the homology of $G$ with arbitrary coefficients.  In
particular, we can use this construction to compute the homology of the
complement of a fiber-type arrangement with coefficients in any local
system. This type of problem has been the focus of a great deal of recent
activity.  In \cite{16}, Esnault, Schechtman, and Viehweg present an
algorithm for computing the cohomology of the complement of an
(arbitrary) arrangement with coefficients in certain complex local
systems (see also \cite{47}).  Refinements of the results of \cite{16}
may be found in \cite{20}, \cite{46}, and \cite{50}.  Note however that
none of these results hold for arbitrary local systems.  See
\cite{24}, \cite{28}, \cite{45}, and~\cite{10} for other results along
similar lines.

The construction of the chain complex $C_{\bullet}(G)$ has certain
functorial properties that allow us to define representations of
groups acting compatibly on an iterated semidirect product of free
groups $G$.  The resulting representations generalize the classical
Magnus representations, \cite{36}, \cite{5}.
Given a compatible automorphism $\psi\in\Aut^\rtimes (G)$, we
explicitly construct the chain equivalence
$\Psi_{\bullet}:C_{\bullet}(G)\to C_{\bullet}(G)$ by means of
``higher-order'' Fox Jacobians of $\psi$. Let $\Gamma$ be a group which
acts on $G$ by compatible automorphisms. Such an action $\Phi:\Gamma\to
\Aut^{\rtimes}(G)$ gives rise to a map $\Phi_{\bullet}$
from $\Gamma$ to the the group of chain automorphisms
of $C_{\bullet}(G)$.  However, this map need not be a homomorphism.
The failure is measured precisely by the chain rule for higher-order
Jacobians (Proposition~4.8).   This problem can be overcome by
following Magnus' original idea in the case $G=F_n$.
Let $\tau:G\to K$ be a $\Phi$-invariant homomorphism.
Then extension of scalars via the map induced by $\tau$ on group rings
yields the desired homomorphism
$\Phi^{\tau}_k: \Gamma\to \Aut_{\Z K}(\Z K \otimes_{\Z G} C_k(G))$, for
each $k\ge 1$ (Theorem 4.11).

In the case where $\Gamma=B_\l$ is the full braid group, acting on
$G=P_{n,\l}$ in a natural fashion, certain choices of homomorphisms
$\tau:P_{n,\l}\to\Z$ yield representations over the ring $\Lambda=\Z\Z$
that generalize the classical Burau representations.   Other
homomorphisms $\tau:P_{n,\l}\to\Z^m$ yield representations of
$B_\l$ that depend on $m$ parameters.  We also obtain generalized
Gassner representations of $P_{\l}$ from the natural action of $P_{\l}$
on $P_{n,\l}$.

Linear representations of braid groups have been much studied over
the years, starting with the pioneering work of Burau, Magnus, and Gassner
(see Birman \cite{5}), and more recently by Jones~\cite{26}, Kohno~\cite{29},
Lawrence~\cite{31}, Moody~\cite{40}, L\"udde and Toppan~\cite{35},
Long and Paton~\cite{34},  Birman, Long and Moody~\cite{6},
and others.  The abiding interest in the subject owes a great
deal to the strong relationship it has with the theory of knots
and links in $S^3$, see \cite{5}, \cite{26}.  The representations
of $B_{\l}$ we obtain here are powerful enough to detect braids
in the kernel of the Burau representation, which is now known
to be unfaithful for $\l\ge 6$, see~\cite{40}, \cite{34}, \cite{6}.
Unlike the representations considered in \cite{26} and \cite{31},
our generalized Burau representations do not factor through the
Hecke algebra in general.  Analogous behavior is exhibited by the
representations constructed in~\cite{29}, \cite{35} and \cite{6}
by other means.

Braid groups are also important in the study of plane algebraic
curves and hyperplane arrangements.  To a curve in
$\C^2$, Moishezon \cite{39} associates a certain ``braid monodromy"
$\theta: F_k\to B_d$, where $d$ is the degree of the curve and $k$
depends on a choice of linear projection $\C^2\to\C$.  Given a
representation $\rho: B_d\to \GL(n,\Lambda)$,
Libgober \cite{33} shows that the $\Lambda$-module
$H_0(F_k; \Lambda^n_{\rho\circ\theta})$ is an invariant of the curve.
In particular, using the reduced Burau representation, he recovers
the Alexander polynomial (up to a factor).  We expect that using the
generalized Burau representations defined here will lead to invariants
of plane curves that cannot be explained solely in terms of the homology
of the maximal abelian cover of the complement.  For an arbitrary
(complexified) hyperplane arrangement, there is an analogous ``pure
braid monodromy,'' see \cite{14}.  We also expect that the
generalized Gassner representations defined here will yield new
invariants of arrangements.

We present several other applications of our construction.  For example,
using it, we obtain algorithms for computing the integral homology
of the Milnor fiber of an arbitrary fiber-type hyperplane
arrangement, as well as the homology eigenspaces of the algebraic
monodromy.  We exhibit the results of some of these computations
in section~7.  The chain complexes arising in these instances
are more manageable than those generated by Aleksandrov~\cite{1} and
Dimca~\cite{15} for the same purposes.  The complexes (of
differential forms) found in these works are generally infinitely
generated.  Also the (rank one) complex local systems arising in the
computation of the eigenspaces of the monodromy are often among those
excluded by the conditions placed on the local system in~\cite{16}.

The homology of the complement of an arrangement with coefficients
in a rank one complex local system is intimately related to the
study of generalized hypergeometric functions,~\cite{2}, \cite{49}.
In light of the work of Schechtman and Varchenko \cite{47},
affine ``discriminantal'' arrangements obtained
from the braid arrangements are of particular interest.  The
fundamental groups of these discriminantal arrangements are of the
form $G=P_{n,\l}$, and may therefore be realized as
iterated semidirect products of free groups.  In section~6, we prove
a vanishing theorem pertaining to these groups which generalizes a
result of Kohno~\cite{29}.  We show that if
$\nu:P_{n,\l}\to\Aut(V)$ is a complex representation (of
arbitrary rank) that is ``quasi-generic'' through rank $q$, then
$H_i(P_{n,\l};V)=0$ for $0\le i\le \min\{q,n-\l-1\}$ (Theorem~6.10).
Results of this form have implications in the study
of the Knizhnik-Zamolodchikov equations, see \cite{47}, \cite{16}.

\remark{Conventions}  Unless otherwise specified, we will regard all
modules over the group ring $\Z G$ of a group $G$ as {\it left} modules.
Elements of the free module $(\Z G)^n$ are viewed as {\it row}
vectors, and $\Z G$-linear maps $(\Z G)^n\to (\Z G)^m$ are
viewed as $n\times m$ matrices which act on the {\it right}
(so that the matrix of $B\circ A$ is $A\cdot B$).  We will
write $[A]^k$ for the map $\oplus_1^k A$ (or,
the block-diagonal $kn\times km$ matrix with diagonal blocks $A$),
$A^{\top}$ for the transpose of $A$, and $I_n$ for the
$n \times n$ identity matrix.

If $U$ and $V$ are two $\Z G$-modules,
$U\otimes_G V$ denotes the $\Z G$-module equal to
$U\otimes V$ modulo the diagonal $G$-action.
If $\phi:G\to H$ is a homomorphism,
$\tilde{\phi}:\Z G\to \Z H$ denotes its extension to group rings,
given by $\tilde{\phi}(\sum n_g g)=\sum n_g\phi(g)$.
(We will abuse notation and also write
$\tilde{\phi}:(\Z G)^n\to (\Z H)^n$ for the map
$\oplus_1^n \tilde\phi$.)
For a $\Z G$-module $V$, there is a $\Z H$-module
$\Z H\otimes_{\Z G} V$ obtained by extension of scalars.
This is achieved by imposing on $\Z H$ the structure of a
{\it right} $\Z G$-module via $s\cdot r=s\tilde\phi(r)$,
and setting $s\cdot(s'\otimes m)=ss'\otimes m$.
An excellent reference for all this is Brown's book \cite{9}.
\endremark

\head 1. Semidirect products of free groups
\endhead

In this section, we introduce the class of groups under consideration
(iterated semidirect products of finitely generated free groups),
give some topological and geometric interpretations,
and provide some examples.

\subhead{1.1}\endsubhead Let $G_1$ and $G_2$ be two groups, and
let $\alpha$ be an action of $G_1$ on $G_2$, i.e., a homomorphism
$\alpha: G_1\to\Aut(G_2)$ from $G_1$ to the group of {\it right}
automorphisms of $G_2$.  The {\it semidirect product} of $G_1$
and $G_2$ with respect to $\alpha$, $G_2\rtimes_{\alpha} G_1$,
is the set $G_2\times G_1$, endowed with the group operation
$(g_2,g_1)\cdot (g'_2,g'_1) = (\alpha(g'_1)(g_2)g'_2,g_1g'_1)$.
The group $G=G_2\rtimes_{\alpha} G_1$ fits into a split
exact sequence
$$1\to G_2 @>\iota_2>> G \overset\pi\to{\underset\iota_1
\to\rightleftarrows} G_1 \to 1,$$
where $\iota_2(g_2)=(g_2,1)$, $\iota_1(g_1)=(1,g_1)$, and
$\pi(g_2,g_1)=g_1$.
Identifying the groups $G_k$ with their images in $G$ under $\iota_k$,
we see that $G$ is generated by $G_1$ and $G_2$, and the following
relations hold in $G$: $g_1^{-1}g_2g_1=\alpha(g_1)(g_2)$, for every
$g_1\in G_1, g_2\in G_2$.

This construction can of course be iterated.  Assume we are given
groups $G_1,\dots, G_{\l}$, and, for each $i<j$, homomorphisms
$\alpha_j^i:G_i \to \Aut(G_j)$ satisfying the compatibility
conditions $\alpha_k^i(g_i)^{-1}\alpha_k^j(g_j)\alpha_k^i(g_i)=
\alpha_k^j(\alpha_k^i(g_i)(g_j))$, for each $i<j<k$.  Then,
we define the {\it iterated semidirect product} of $G_1,\dots, G_{\l}$
with respect to the actions $\alpha_j^i$ to be the group
$$G=G_\l \rtimes_{\alpha_{\l}} G_{\l-1} \rtimes
\dots \rtimes_{\alpha_{3}} G_2 \rtimes_{\alpha_{2}} G_1,$$
where, for each $1\le q \le \l$, the partial iteration,
$G^q=G_q\rtimes_{\alpha_q} G^{q-1}$, is defined by the homomorphism
$\alpha_q: G^{q-1}\to \Aut(G_q)$, whose restriction to $G_p$,
$1\le p<q$, is $\alpha_q^p$.

In this paper, we study in detail groups $G$ which may be realized as
iterated semidirect products of finitely generated free groups.
Such groups can be written as $G = \rtimes_{q=1}^{\l} G_q$, where
$G_q = F_{d_q} = \langle x_{1,q},\dots,x_{d_q,q}\rangle$
is free on $d_q$ generators.
It follows readily that the group $G$ has presentation
$$G=\langle x_{i,q}\quad (1\le i \le d_q, 1\le q \le \l) \mid
x_{j,p}^{-1}x_{i,q}x_{j,p}=\alpha_q^{j,p}(x_{i,q})\quad
(p<q) \rangle,\tag{1.2}$$
where $\alpha_q^{j,p}:=\alpha_q(x_{j,p})\in \Aut(F_{d_q})$.
Conversely, any group $G$ with presentation as above admits
the structure of an iterated semidirect product of free groups
in an obvious fashion.

\example{Example 1.3} The principal motivation for our
analysis of iterated semidirect products of free groups are
Artin's (pure) braid groups.  Let
$$B_\l=\langle \sigma_i\ \ (1\le i<\l) \mid
\sigma_i\sigma_{i+1}\sigma_i=\sigma_{i+1}\sigma_i\sigma_{i+1}
\ (1\le i<\l-1),\ \sigma_i\sigma_j=\sigma_j\sigma_i
\ (|i-j|>1)\rangle$$ denote the braid
group on $\l$ strings, and $P_\l$ the subgroup of braids
with trivial permutation of the strings, see Artin~\cite{4},
and the books by Birman~\cite{5} and Hansen~\cite{23}.
The pure braid group,
$P_{\l}=F_{\l-1}\rtimes_{\alpha_{\l-1}}\cdots\rtimes_{\alpha_2} F_1$,
admits the structure of an
iterated semidirect product of free groups.
The monodromy homomorphisms
$\alpha_q:P_{q}\to \Aut(F_q)$, $2\le q\le \l-1$ are
given by the restriction to $P_q$ of the Artin representation,
$\alpha_q: B_{q}\to \Aut(F_q)$, defined by
$$\alpha_q(\sigma_i)(x_j)=\cases
x_i x_{i+1} x_i^{-1}&\text{if $j=i$,}\\
x_i&\text{if $j=i+1$,}\\
x_j&\text{otherwise,}
\endcases$$
where $F_q=\langle x_1,\dots,x_q\rangle$.
The iterated semidirect product structure is in evidence in the familiar
presentation of the pure braid group found in the above references. The
group $P_{\l}$ has generators $$A_{i,j}=\sigma_{j-1}\cdots\sigma_{i+1}
\sigma_i^2\sigma_{i+1}^{-1}\cdots\sigma_{j-1}^{-1},\quad
1\le i < j \le \l,$$ and defining relations
$$A_{r,s}^{-1} A_{i,j} A_{r,s} =
\cases A_{i,j}&\text{if $i<r<s<j$ or $r<s<i<j$,}\\
A_{r,j}A_{i,j}A_{r,j}^{-1}&\text{if $r<s=i<j$,}\\
A_{r,j}A_{s,j}A_{i,j}A_{s,j}^{-1}A_{r,j}^{-1}&\text{if $r=i<s<j$,}\\
[A_{r,j},A_{s,j}]A_{i,j}[A_{r,j},A_{s,j}]^{-1}&\text{if $r<i<s<j$.}\\
\endcases$$
\endexample

\subhead{1.4}\endsubhead Let us give a topological interpretation of
iterated semidirect products of (finitely generated) free groups.
Associated to a group $G = \rtimes_{i=1}^{\l} F_{d_i}$ there is a
standard CW-complex $X=X_G$, with fundamental group isomorphic to $G$.
This complex is defined inductively as a tower of fibrations,
$$X=X^{\l} @>p_{\l}>> X^{\l -1} @>p_{\l -1}>> \cdots @>p_3>> X^2
@>p_2>> X^1,$$
such that each projection $p_q:X^q\to X^{q-1}$ admits a section,
and has fiber $K_{d_q}=\bigvee_1^{d_q} S^1$.

The complex $X$ is constructed as follows:  Take $X^1=K_{d_1}$.
Inductively assume that the space
$X^{q-1}$ with $\pi_1(X^{q-1})=G^{q-1}$ has been constructed.
Let $\E_0(K)$ denote the group of based homotopy
classes of based self-homotopy equivalences of a space $(K,y_0)$.
If $K$ is a $K(\pi, 1)$ space, the evaluation map
$\ev: \E_0(K) \to \Aut(\pi_1(K, y_0))$, defined by $\ev(f)=f_{\#}$, is
an isomorphism.  Applying this observation to $K=K_{d_q}$, we
see that the homomorphism $\alpha_q:G^{q-1} \to \Aut(F_{d_q})$
factors through $\tau_q:G^{q-1} \to \E_0(K_{d_q})$.
Now let $\pi:\widetilde{X}^{q-1}\to X^{q-1}$ be the universal cover
of $X^{q-1}$, and identify $G^{q-1}$ as the group of deck
transformations.  Consider the diagonal action of $G^{q-1}$ on
$\widetilde{X}^{q-1}\times K_{d_q}$, defined by
$g\cdot (\tilde{x}, y) = (g\tilde{x}, \tau_q(g)(y))$, and
let $X^{q}$ denote the orbit space of this action.
By construction, the map $p_q:X^{q}\to X^{q-1}$ given by
$p_q([\tilde{x}, y])=\pi(\tilde{x})$ is a fibration, with
fiber $K_{d_q}$.  Moreover, a canonical section
$s_q: X^{q-1}\to X^{q}$ is given by $s_q(x)=[\tilde{x},y_0].$
(This is well-defined, as $[\tilde{x},y_0]=[g\tilde{x},
\tau_q(g)(y_0)]=[g\tilde{x}, y_0]$.)  Finally,
$\pi_1(X^{q}) = F_{d_q} \rtimes_{\alpha_q} G^{q-1}=G^q$,
and this finishes the inductive construction of $X=X^{\l}$.

Notice that the CW-complex $X_G$ defined above is a $K(G,1)$
space.   This follows from the long exact sequence in homotopy
for a fibration and induction.  Since $X_G$ is $\l$-dimensional
by construction, the group $G$ is of type FL, and its cohomological
dimension is $\l$, see \cite{9}.

\subhead{1.5}\endsubhead A more geometric interpretation of the
group $G$ is when the above tower of (Serre) fibrations can be
replaced (up to homotopy) by a tower of locally trivial bundles,
$$M=M^{\l} @>\pi_{\l}>> M^{\l -1} @>\pi_{\l -1}>> \cdots
@>\pi_3>> M^2 @>\pi_2>> M^1,$$
such that the fiber $\Sigma_q$ of $\pi_q:M^q\to M^{q-1}$ is a
surface with a number of punctures. In this case, the homomorphism
$\alpha_q:G^{q-1} \to \Aut(F_{d_q})$ can be realized by the
monodromy of the bundle, $\mu_q:M^{q-1}\to \Homeo(\Sigma_q)$.

Perhaps the simplest situation is where each fiber has genus 0,
i.e.,~$\Sigma_q = \C \setminus \{ d_q \text{ points}\}$.  This is the
case, for example, when $M$ is the complement of a fiber-type
arrangement of complex hyperplanes.  In this instance,
each map $\pi_q:M^q\to M^{q-1}$ is, by definition, the restriction
of a linear projection $\C^{q}\to\C^{q-1}$.
This notion was introduced by Falk and Randell in \cite{18},
where they prove the LCS formula for $G=\pi_1(M)$, the group of a
central fiber-type arrangement.  This result was subsequently extended to
arbitrary fiber-type arrangements by Jambu \cite{25}.
The {\it exponents} $\{d_1,\dots,d_\l\}$ arising in the iterated semidirect
product structure on $G$ from the iterated bundle structure on $M$ (i.e.~the
exponents of the fiber-type arrangement itself) are, via the LCS formula,
determined by the Betti numbers of $G$.
In the case of the braid arrangement, $\A_\l=\{H_{i,j}=\ker(z_i-z_j)\}$,
whose complement, $M(\A_\l)=\C^\l\setminus\bigcup H_{i,j}$,
is the configuration space of the set of
$\l$ (ordered) points in $\C$, the fiber-type structure was
first discovered by Fadell and Neuwirth \cite{17}.  The resulting
decomposition of
$P_{\l}=\pi_1(M(\A_\l))$ is precisely the one exhibited in Example~1.2.

\subhead{1.6}\endsubhead We conclude this section with a few remarks
concerning the non-uniqueness of semidirect product structures of groups.
First, note that if $\alpha,\beta:G_1\to \Aut(G_2)$ are
homomorphisms that differ by an {\it inner} automorphism $\gamma$
of $G_2$ (i.e., $\alpha(g)=\gamma\cdot\beta(g)$,
$\forall g\in G_1$), then the semidirect products
$G_2\rtimes_{\alpha} G_1$ and $G_2\rtimes_{\beta} G_1$ are isomorphic.
Thus, the isomorphism class of the group
$G=G_2\rtimes_{\alpha} G_1$ depends only on the homomorphism
$\bar{\alpha}:G_1\to\operatorname{Out}(G_2)$.
Second, let us point out that there is no well-defined notion of
exponents of a group in general. That is, for a group $G$ that can
be written as $G = \rtimes_{i=1}^{\l} F_{d_i}$, the
``exponents'' $\{ d_1,\dots,d_\l\}$  depend on
the particular iterated semidirect product structure on $G$.
(On the other hand, the number of exponents, $\l$, depends
only on $G$, since $\operatorname{cd}(G)=\l$.)
We illustrate these phenomena with two examples that are
relevant to our general discussion.

\example{Example 1.7} The pure braid group on 3 strings, $P_3$.
It follows from Example~1.3 that $P_3=F_2\rtimes_{\alpha_2} F_1$, where
$F_1=\langle A_{1,2} \rangle$, $F_2=\langle A_{1,3}, A_{2,3}
\rangle$, and $\alpha_2(A_{1,2})$ is conjugation by
$A_{1,3}A_{2,3}$.  Thus, $P_3$ is isomorphic to the direct
product $F_2\times F_1$.

There is another realization of $P_3$ as a
semidirect product of free groups:  $P_3=F_4\rtimes_{\mu} F_1$,
corresponding to the Milnor fibration of the braid arrangement
in $\C^3$ (see section~7).  A computation shows that
$F_1=\langle A_{1,2} \rangle$,
$F_4=\langle t_1, t_2, t_3, t_4 \rangle$, where
$t_1=A_{1,2}^{-1}A_{1,3},\, t_2=A_{1,2}^{-1}A_{2,3},\,
t_3=A_{1,3}A_{1,2}^{-1},\, t_4=A_{2,3}A_{1,2}^{-1}$,
and the action $\mu:F_1\to \Aut(F_4)$ is given by
$$
A_{1,2}:\cases
t_1 \mapsto t_1t_4t_2^{-1}\\
t_2 \mapsto t_1t_4t_3^{-1}\\
t_3 \mapsto t_1\\
t_4 \mapsto t_2
\endcases
$$
\endexample

\example{Example 1.8} The pure braid group on 4 strings, $P_4$.
As discussed in Example~1.3, this group may be
realized as $P_4\cong F_3 \rtimes_{\alpha_3} F_2 \rtimes_{\alpha_2} F_1$.
Since the Coxeter groups $\text{A}_3$ and $\text{D}_3$ are
isomorphic, we have $P_4 = \PB(\text{A}_3) \cong
\PB(\text{D}_3)$.  This latter group may be realized as
$\PB(\text{D}_3)\cong F_5\rtimes_{\beta_3} F_2\rtimes_{\alpha_2} F_1$.

The geometric reason for this decomposition is due to Brieskorn
\cite{8}, who found a (non-linear) bundle map from the complement
of the $\text{D}_\l$ arrangement to a hyperplane complement homotopy
equivalent to the complement of the
$\text{A}_{\l-1}$ arrangement.  This map was studied by Falk and Randell
\cite{18}, who noted that the Brieskorn bundle admits a section, and
that its fiber is a curve of genus $2^{\l-2}(\l-3)+1$ with $2^{\l-1}$
punctures.
It follows that $\PB(\text{D}_\l) \cong F_k \rtimes_{\beta_\l} P_\l$,
where $k=2^{\l-1}(\l-2)+1$.  The representation
$\beta_\l: P_\l\to \Aut(F_k)$
was recently identified by Leibman and Markushevich \cite{32}.
The action of $P_3$, generated by $\{A_{1,2},A_{1,3},A_{2,3}\}$,
on $F_5=\langle t_1,\dots,t_5\rangle$ is given by:
$$
A_{1,2}: \cases
t_1 \mapsto t_1\\
t_2 \mapsto t_2t_4t_5^{-1}t_3^{-1}t_2t_1\\
t_3 \mapsto t_2t_4t_5^{-1}t_1\\
t_4 \mapsto t_1^{-1}t_2^{-1}t_3t_5\\
t_5 \mapsto t_1^{-1}t_5t_4^{-1}t_2^{-1}t_3t_5
\endcases
\qquad
A_{1,3}: \cases
t_1 \mapsto t_2^{-1}t_5\\
t_2 \mapsto t_2\\
t_3 \mapsto t_3t_1t_5^{-1}t_2\\
t_4 \mapsto t_2^{-1}t_4t_1^{-1}t_5\\
t_5 \mapsto t_2^{-1}t_5t_1^{-1}t_5
\endcases
\qquad
A_{2,3}: \cases
t_1 \mapsto t_3^{-1}t_1t_4\\
t_2 \mapsto t_3^{-1}t_2t_4\\
t_3 \mapsto t_3\\
t_4 \mapsto t_4\\
t_5 \mapsto t_5
\endcases
$$
(This is not how the representation $\beta_3:P_3\to \Aut(F_5)$ is
written in \cite{32}; their formula for $A_{1,3}$ is not correct, but
can be fixed by carefully following their algorithm.)
\endexample

\head 2. The Resolution
\endhead

In this section, we construct a finite, free $\Z G$-resolution
of the integers for every group $G$ which admits the structure
of an iterated semidirect product of finitely generated free groups.
The basis for this construction is the non-commutative differential
calculus for words in a free group developed by Fox in \cite{21}
(see \cite{5} for an exposition).

\subhead{2.1}\endsubhead First consider a single free group
$F_n=\langle x_1,\dots,x_n \rangle$.
Let $K_n=\bigvee_1^{n} S^1$ be the standard $K(F_n,1)$, and let
$\widetilde C_{\bullet}$ be the augmented chain complex of
the universal cover $\widetilde K_n$.  Identifying $C_0$
with $\Z F_n$, and $C_1$ with $(\Z F_n)^n$ (with basis
$\{e_1,\dots , e_n\}$ given by the lifts of the 1-cells at the
basepoint), the resolution $\widetilde C_{\bullet}$ can
be written as:
$$0\to (\Z F_n)^n @>\Delta>> \Z F_n @>\epsilon>> \Z \to 0,$$
where $\Delta=\pmatrix x_1-1 &\cdots & x_n-1\endpmatrix^{\top}$
and $\epsilon$ is the augmentation map, given by $\epsilon(x_i)=1$.
Further,
consider an automorphism $\alpha:F_n\to F_n$.  The induced chain map
$\alpha_{\bullet}:C_{\bullet} \to C_{\bullet}$ can be written as:
$$\CD
(\Z F_n)^n        @>\Delta>>             \Z F_n\\
@VVJ(\alpha)\circ \tilde\alpha V   @VV\tilde\alpha V \\
(\Z F_n)^n        @>\Delta>>             \Z F_n
\endCD\tag{2.2}
$$
where $\DS{J(\alpha) =
\fracwithdelims(){\partial\alpha(x_i)}{\partial x_j}}$
is the $n\times n$ Jacobian matrix of Fox derivatives of $\alpha$.
Note that $\alpha_{\bullet}$ is the composition
of a $\Z F_n$-linear map ($\id_{\Z F_n}$, resp. $J(\alpha)$),
with a non-linear map (the extension $\tilde\alpha:\Z F_n\to \Z F_n$,
resp. $\tilde\alpha:(\Z F_n)^n\to (\Z F_n)^n$).
The commutativity of diagram \thetag{2.2} is a consequence of
the ``fundamental formula of Fox Calculus.''

If $\beta:F_n\to F_n$ is another automorphism, the fact that
$(\beta\circ\alpha)_{\bullet}=\beta_{\bullet}\circ\alpha_{\bullet}$
is a consequence of the following ``chain rule of Fox Calculus:''
\proclaim{Lemma 2.3}
$J(\beta\circ\alpha)=\tilde\beta(J(\alpha))\cdot J(\beta)$.
\endproclaim
\noindent
In particular, $J(\alpha^{-1})=\tilde\alpha^{-1}(J(\alpha)^{-1})$.

\subhead{2.4}\endsubhead Given $F_n=\langle x_1,\dots, x_n\rangle$ and
$\alpha\in\Aut(F_n)$ as above, form the semidirect product,
$G_\alpha := F_n \rtimes_{\alpha} F_1=
\langle x_i,t \mid t^{-1}x_it=\alpha(x_i) \rangle$, of $F_1=\langle t \rangle$
with $F_n$ determined by $\alpha$.
Let $R=\Z G_\alpha$, and define $\lambda_t: R\to R$ by
$\lambda_t(g)=t\cdot g$.
(Note that $\lambda_t$ is {\it not} $R$-linear with respect to the
left-module structure on $R$, unless $G_{\alpha}$ is abelian.)
Extension of scalars yields the following commuting diagram:
$$\CD
R\otimes_{\Z F_n} (\Z F_n)^n @>{\id\otimes\Delta}>>
R\otimes_{\Z F_n} \Z F_n \\
@VV\lambda_t \otimes J(\alpha)\circ\tilde\alpha V
@VV\lambda_t \otimes \tilde\alpha V \\
R\otimes_{\Z F_n} (\Z F_n)^n @>{\id\otimes\Delta}>>
R\otimes_{\Z F_n} \Z F_n \\
\endCD
$$
The map $\lambda_t \otimes J(\alpha)\circ\tilde\alpha$, together with the
canonical isomorphism $R\otimes_{\Z F_n} (\Z F_n)^n \cong R^n$,
define a map $\rho(\alpha):R^n \to R^n$.

\proclaim{Lemma 2.5}  The map $\rho(\alpha)$ belongs to
$\Aut_R(R^n)$, and its matrix is $t\cdot J(\alpha)$.
\endproclaim

\demo{Proof}  First we must verify that $\rho(\alpha)$
is $R$-linear.  For that, start by computing $\rho(\alpha)$
on the basis vectors $\{e_1,\dots , e_n\}$:
$$\rho(\alpha)(e_i)=t \sum_{j=1}^n \frac{\partial
\alpha(x_i)}{\partial x_j} \cdot e_j.\tag{$\dagger$}$$
Clearly, $\rho(\alpha)(t^k\cdot e_i)= t^k\cdot \rho(\alpha)(e_i)$.
For an element $w$ of $F_n$, we have
$$
\rho(\alpha)(w\cdot e_i)=t\cdot\alpha(w)  \sum_{j=1}^n
\frac{\partial \alpha(x_i)}{\partial x_j} \cdot e_j
=t\cdot t^{-1}wt \sum_{j=1}^n \frac{\partial
\alpha(x_i)}{\partial x_j} \cdot e_j
=w\cdot \rho(\alpha)(e_i),
$$
and $R$-linearity follows from the fact that $F_1$ and $F_n$
generate $G_\alpha=F_n\rtimes_{\alpha} F_1$.  That the matrix of
$\rho(\alpha)$ is as asserted follows from \thetag{$\dagger$}.

Finally, we must verify that $\rho(\alpha)$ has an inverse.
Note that $G_\alpha\cong G_{\alpha^{-1}}$, the isomorphism
being given by $x_i\mapsto x_i$ and $t\mapsto t^{-1}$.
Thus, $\Z G_{\alpha^{-1}} = R$.  It now follows from Lemma 2.3
that $\rho(\alpha)\circ\rho(\alpha^{-1}) = \id$.
\quad\qed\enddemo

\subhead{2.6}\endsubhead We now consider an iterated semidirect product
$G \cong G_\l \rtimes_{\alpha_{\l}}
\dots \rtimes_{\alpha_3} G_2\rtimes_{\alpha_2} G_1$, where
$G_q = F_{d_q} = \langle x_{1,q},\dots,x_{d_q,q}\rangle.$  Recall that
$G^q$ denotes the split extension $G_q\rtimes_{\alpha_q} G^{q-1}$,
and $\alpha^p_q: G_p\to \Aut(G_q)$ denotes the restriction of
$\alpha_q$ to $G_p$, $1\le p<q$.
Let $R=\Z G$ denote the integral group ring of $G$.
For each generator $x_{r,p}$ $(1\le p<q, 1\le r\le d_p)$ of
$G^{q-1}$, we have a commuting diagram
$$\CD
(\Z G_q)^{d_q} @>{\Delta_q}>> \Z G_q \\
@VV{J^{r,p}_q}\circ\tilde\alpha_q^{r,p}V  @VV{\tilde\alpha_q^{r,p}}V \\
(\Z G_q)^{d_q} @>{\Delta_q}>> \Z G_q
\endCD
$$
where
$\Delta_q=\pmatrix x_{1,q}-1 &\cdots & x_{d_q,q}-1\endpmatrix^{\top}$,
$\tilde\alpha_q^{r,p}$ is the extension of the automorphism
$\alpha_q^{r,p}:=\alpha_q^p(x_{r,p}):G_q\to G_q$ induced by conjugation
by $x_{r,p}$, and $J^{r,p}_q = J(\alpha_q^{r,p})$ is the
Jacobian matrix of $\alpha_q^{r,p}$.
Let $\lambda_{r,p}:R\to R$ be left multiplication by $x_{r,p}$.
By Lemma 2.5, the map
$$\rho(\alpha_q^{r,p})=\lambda_{r,p}\otimes
J^{r,p}_q\circ\tilde\alpha_q^{r,p}: R\otimes_{\Z G_q} (\Z G_q)^{d_q}
\to R\otimes_{\Z G_q} (\Z G_q)^{d_q},$$
defines an $R$-linear automorphism of $R^{d_q}$, with matrix
$\DS{x_{r,p}\cdot\left({\partial
\alpha_q^{r,p}(x_{i,q}) \over \partial x_{j,q}}\right).}$

\proclaim{Lemma 2.7} For each $q$, $1<q\le \l$, there is a
(unique) representation $\rho_q:G^{q-1}\to \Aut_R(R^{d_q})$ with
the property that $\rho_q(x)=\lambda_x\otimes
J(\alpha_q(x))\circ\tilde\alpha_q(x)$ for every $x \in G^{q-1}$.
\endproclaim
\demo{Proof}   Specifying the automorphisms $\rho_q(x_{r,p})$
for each $1\le r\le d_p$ defines a
representation $\rho^p_q:G_p\to \Aut_R(R^{d_q})$ for each of the
free group $G_p=F_{d_p}$, $1\le p < q\le \l$.  We are left with
showing that these representations are compatible with the
iterated semidirect product structure on $G$, i.e.,
$\rho^p_s(x_{j,p}^{-1})\rho^q_s(x_{i,q})\rho^p_s(x_{j,p})=
\rho^q_s(\alpha_q^{j,p}(x_{i,q}))$, for $p<q<s$.  This
follows from the fact that $\lambda_x\circ \lambda_y=\lambda_{xy}$
and Lemma 2.3.
\quad\qed\enddemo

\subhead 2.8\endsubhead  The above representation extends to a
representation $\rho_q:G\to \Aut_R(R^{d_q})$ via the convention
$\rho_q(x_{r,p}) = I_{d_q}$ if $p\ge q$.   We denote by
$\tilde\rho_q:R\to \End_R(R^{d_q})$ the extension of
$\rho_q$ to the group ring $R$.  Replacing each entry $x$ of an
$m \times n$ matrix by $\tilde\rho_q(x)$ defines a homomorphism
$\Hom_R(R^m,R^n)\to\Hom_R(R^{md_q},R^{nd_q})$ that we still
denote by $\tilde\rho_q$.  By restriction, we also get a
homomorphism $\tilde\rho_q:\Aut_R(R^n)\to \Aut_R(R^{nd_q})$.

\subhead{2.9}\endsubhead  We now construct a free resolution
$\epsilon:C_\bullet\to \Z$ over $R=\Z G$.  Let $C_0=R$, and for
$1\le k\le\l$ let
$$C_k = \bigoplus_{1\le p_1 <\dots<p_k\le\l}
R^{d_{p_1}d_{p_2}\cdots d_{p_k}}.$$
The augmentation map, $\epsilon: C_0\to \Z$, is the
usual augmentation of the group ring, given by $\epsilon(g)=1$,
for $g\in G$.  We define the boundary maps of
the complex $C_{\bullet}$ by specifying their restrictions
$\Delta^{p_1,p_2,\dots,p_k}$ to the summands
$R^{d_{p_1}d_{p_2}\cdots d_{p_k}}$.  This is done recursively as
follows:
\smallskip
\noindent Define $\Delta_p:R^{d_p}\to R$ by $\Delta_p =
\pmatrix x_{1,p}-1 &\cdots & x_{d_p,p}-1 \endpmatrix^{\top}$.
\smallskip
\noindent For $p_1<p_2$, define
$\Delta_{p_1,p_2}:R^{d_{p_1}d_{p_2}}\to R^{d_{p_2}}$ by
$\Delta_{p_1,p_2} =-\tilde\rho_{p_2}(\Delta_{p_1})$.
\smallskip
\noindent In general, for $1\le p_1 <\dots<p_k\le\l$, define
$\DS{\Delta_{p_1,\dots,p_k}:R^{d_{p_1}\cdots d_{p_k}}\to
R^{d_{p_2}\cdots d_{p_k}}}$ by $$\Delta_{p_1,\dots,p_k} =
-\tilde\rho_{p_k}(\Delta_{p_1\dots,p_{k-1}}).$$
\smallskip
\noindent Now define $\Delta^{p_1,\dots,p_k}: R^{d_{p_1}\cdots
d_{p_k}}\to \bigoplus_{i=1}^k  R^{d_{p_1}\cdots \hat d_{p_i}\cdots
d_{p_k}}$  by $$\Delta^{p_1,\dots,p_k} =
\left(\Delta_{p_1\dots,p_{k}},
\left[\Delta_{p_2,\dots,p_k}\right]^{d_{p_1}},\dots,
\left[\Delta_{p_i,\dots,p_k}\right]^{d_{p_1}\cdots d_{p_{i-1}}},
\dots, \left[\Delta_{p_k}\right]^{d_{p_1}\cdots d_{p_{k-1}}}\right).$$
\smallskip
\noindent Finally, define $\Delta:C_k\to C_{k-1}$ by
$\DS{\Delta = \bigoplus_{1\le p_1 <\dots<p_k\le\l}
\Delta^{p_1,\dots,p_k}}$.

In the context of the above construction, the fundamental formula of
Fox calculus has the following consequences.

\proclaim{Lemma 2.10} For $x\in G_i$ and
$1\le i<p_1<\cdots<p_k<q\le\l$, we have
$$\tilde\rho_{q}\circ\tilde\rho_{p_{k}}\circ\cdots\circ\tilde\rho_{p_1}(x)
\cdot [\Delta_q]^{d_{p_1}\cdots d_{p_k}} = [\Delta_q]^{d_{p_1}\cdots
d_{p_k}}\cdot
\tilde\rho_{p_{k}}\circ\cdots\circ\tilde\rho_{p_1}(x).$$
\endproclaim
\demo{Proof} First consider the case $k=0$.  Let $\alpha\in\Aut(G_q)$
be the automorphism induced by conjugation by $x$.  Then the matrix of
$\rho_q(x)$ is $x\cdot J(\alpha)$ (see~2.6).  Using the fundamental
formula of Fox calculus as in 2.2, we have
$$\rho_q(x)\cdot\Delta_q=x\cdot J(\alpha)\cdot\Delta_q =
x\cdot\tilde\alpha(\Delta_q) = \Delta_q\cdot x.$$

In general, write
$A=\tilde\rho_{p_{k}}\circ\cdots\circ\tilde\rho_{p_1}(x)$ and note
that $A$ is a square matrix of size $d=d_{p_1}\cdots d_{p_k}$.  For each
entry $a$ of $A$, we have $\tilde\rho_q(a)\cdot\Delta_q=
\Delta_q\cdot a$ by 2.2 as above.  It then follows from some
elementary matrix manipulations that $\tilde\rho_q(A)\cdot[\Delta_q]^d
= [\Delta_q]^d \cdot A$.
\quad\qed
\enddemo

We now come to the main theorem of this section.

\proclaim{Theorem 2.11}  Given a group $G$ which admits the
structure of an iterated semidirect product of finitely
generated free groups, the system of $R$-modules and homomorphisms
$\{C_\bullet,\Delta\}$ is a finite, free resolution of $\Z$ over
$R=\Z G$.
\endproclaim
\demo{Proof}  The proof is by induction on $\l$ with the case
$\l =1$ clear.

Let $G=\rtimes_{p=1}^{\l} G_p$, where $G_p = F_{d_p}$, and consider
the  (normal) subgroup $\Cal G<G$ given by $\Cal G=\rtimes_{p=2}^{\l}
G_p$.  By induction, the construction of 2.9 yields a free
resolution $\varepsilon: C_\bullet(\Cal G)\to\Z$ over $\Cal R=\Z \Cal G$.

Let $\widehat C_\bullet= C_\bullet(\Cal G)\otimes_{\Cal G} R$, and let
$D_\bullet$ denote the chain complex of $R$-modules with terms $D_k =
( \widehat C)^{d_1}$ and differentials $\partial_D =
-[\partial_{\widehat C}]^{d_1}$. That is, $D_\bullet$ is the direct
sum of $d_1$ copies of $\widehat C_\bullet$, with the sign of the
differential reversed.  Note that
$\widehat C_\bullet$ and $D_\bullet$ are complexes of free
$R$-modules, and that $H_*(\Cal G;R) = H_*(\widehat C_\bullet)$.

Define a map
$\Xi_\bullet:D_\bullet\to\widehat C_\bullet$ by setting the restriction
of $\Xi_\bullet$ to the summand $R^{d_1d_{p_1}\cdots d_{p_k}}$ of
$D_k$ to be equal to
$$\Delta_{1,p_1,\dots,p_k}:R^{d_1d_{p_1}\cdots d_{p_k}}\to
R^{d_{p_1}\cdots d_{p_k}}\subset\widehat C_k.$$ In particular,
$\Xi_0:D_0\to\widehat C_0$ is given by
$\Xi_0=\Delta_1:R^{d_1}\to R$.  We claim that
$\Xi_\bullet:D_\bullet\to\widehat C_\bullet$ is a chain map. To prove
this assertion, it suffices to verify that
$$\Delta_{1,p_1,\dots,p_k}\cdot [\Delta_{p_j,\dots,p_k}]^{d_{p_1}\cdots
d_{p_{j-1}}} =  -[\Delta_{p_j,\dots,p_k}]^{d_1d_{p_1}\cdots
d_{p_{j-1}}}\cdot\Delta_{1,p_1,\dots,\hat p_j,\dots,p_k}$$
for $1\le j\le k$ (where $d_{p_1}\cdots d_{p_{j-1}}=1$ if $j=1$).
These equalities all follow easily from Lemma~2.10.
Furthermore, it is clear from the construction in 2.9 that
$C_\bullet=C_\bullet(G)$ is the mapping cone of the chain map
$\Xi_\bullet$.  Thus $C_\bullet$ is a chain complex of free
$R$-modules.

We now show that $C_\bullet\to\Z$ is a resolution, i.e.~that
$\widetilde C_\bullet$ is acyclic.  It is clear from the construction
that $H_0(\widetilde C_\bullet)=0$.  Since $C_\bullet$ is the mapping
cone of $\Xi_\bullet:D_\bullet\to\widehat C_\bullet$, we have a
long exact sequence in homology
$$\cdots\to H_{i+1}(C_\bullet)\to H_i(D_\bullet)@>{H_i(\Xi_\bullet})>>
H_i(\widehat C_\bullet)\to  H_i(C_\bullet)\to\cdots$$ with connecting
homomorphisms induced by the chain map $\Xi_\bullet$. Now $R$ is free
as a $\Cal G$-module and $H_*(\Cal G;R) = H_*(\widehat C_\bullet)$, so
$$H_i(\widehat C_\bullet) =
\cases R_{\Cal G} &\text{if $i=0$,} \\  0 &\text{if $i\ne 0$,}
\endcases
\qquad\text{and therefore}\qquad H_i(D_\bullet)=\cases (R_{\Cal
G})^{d_1} &\text{if $i=0$,} \\ 0 &\text{if $i\ne 0$.}
\endcases
$$ Thus the long exact sequence above reduces to
$$0\to H_1(C_\bullet)\to (R_{\Cal G})^{d_1} @>{H_0(\Xi_\bullet)}>>
R_{\Cal G}\to H_0(C_\bullet)\to 0,$$ and we have $H_i(C_\bullet)=0$ for
$i\ge 2$.  It remains to show that $H_1(C_\bullet)=0$, i.e.~that the
connecting homomorphism $H_0(\Xi_\bullet)$ is injective.

Now $H_0(\widehat C_\bullet) = R_{\Cal G} =
\Z G\otimes_{ \Cal G}\Z = \Z [G/\Cal G] = \Z F_{d_1}$, so
$H_0(D_\bullet)=(\Z F_{d_1})^{d_1}$. Under this identification we have
$H_0(\Xi_\bullet) = (x_{1,1}-1\ \cdots\ x_{d_1,1}-1)^{\top}$, where
$F_{d_1} =
\langle x_{1,1},\dots,x_{d_1,1}\rangle$. Hence
$H_0(\Xi_\bullet)$ is injective, and the proof is complete.
\quad\qed

\enddemo

\remark{Remark 2.12} The chain complex $C_{\bullet}(G)$ is,
by acyclic models, chain-equivalent to $C_{\bullet}(\widetilde{X}_G)$,
the equivariant chain complex of the universal cover of the $K(G,1)$ space
$X_G$ constructed in 1.4.  For example, if $G = F_n\rtimes_{\alpha} F_m$
is the semidirect product of {\it two} free groups, it is immediate from
the construction of 2.9 that $C_\bullet(G)$ is the chain
complex resulting from application of the Fox Calculus to the
presentation
$G = \langle x_i,y_j\mid y_j x_i = x_i \alpha(x_i)(y_j)\rangle$,
where $F_m=\langle x_i\rangle$ and $F_n=\langle y_j\rangle$.
Thus, in this instance, $C_\bullet(G)$ is equal to the chain complex of the
universal cover of the ``presentation two-complex'' $X_G$, associated to
this presentation of the group $G$.
\endremark

\remark{Remark 2.13}  After completing this work, we became aware
of a construction of Brady~\cite{7}.  Given a semidirect product
$G=G_2\rtimes_{\alpha} G_1$, a free resolution of $G_1$, and a
free resolution of $G_2$ which admits an action of $G_1$ compatible
with $\alpha$, Brady describes in \cite{7} an algorithm for producing a
free resolution of $G$.  Carrying out this algorithm inductively
in our situation, and making use of the lemmas from 2.1--2.8, it is
possible to show that the resulting resolution coincides with the one
constructed in 2.9.  Although this argument is somewhat
shorter than the one presented here, the explicit construction and
arguments presented above will be of further use in the remainder of
the paper.
\endremark

\head 3. Some Consequences
\endhead

In this section we derive some consequences of Theorem~2.11,
and illustrate the construction of the previous section by
means of several examples.

\subhead 3.1. Direct Products \endsubhead  We first consider the simplest
situation---that of a direct product.
\proclaim{Proposition 3.2}  Let $G$ and $H$ be two iterated
semidirect products of free groups, and let $C_\bullet(G)$ and
$C_\bullet(H)$ be the corresponding chain complexes.  Then
$G\times H$ is also an iterated semidirect product of free
groups, and
$C_\bullet(G\times H)=C_\bullet(G)\otimes C_\bullet(H).$
\endproclaim

\demo{Proof}  The iterated semidirect product structure on
$G\times H$ is obtained in an obvious fashion.  The structure
of $C_\bullet(G\times H)$ follows from the construction
and standard facts about tensor products of resolutions (see
e.g. \cite{9}, p. 107.)\quad\qed
\enddemo

This result has the following topological interpretation.
Recall the  $K(G,1)$ space $X_G$ defined in section~1.
It follows from the construction that
$X_{G\times H} \simeq X_G \times X_H$.  Passing to
equivariant chain complexes of universal covers recovers
the above result.

\example{Example 3.3} Let $G=\times_1^\l G_p$, where $G_p=F_{d_p}$.
Then $\widetilde C_{\bullet}(G)= \otimes_1^\l
\widetilde C_{\bullet}(G_p)$, where $\widetilde C_{\bullet}(G_p):
(\Z G_p)^{d_p} @>\Delta^{(p)}>>\Z G_p @>\epsilon^{(p)}>>\Z\to 0$
is as in 2.1.  Explicitly, $C_k(G)$ is the direct sum of
$C_{i_1}(G_1)\otimes\cdots \otimes  C_{i_\l}(G_\l)$,
over all indices $i_r\in\{0,1\}$ such that $i_1+\cdots + i_\l=k$,
and the restriction of the differential $\Delta:C_k(G)\to C_{k-1}(G)$
to such a summand is given by
$$\Delta(c_1\otimes\cdots\otimes c_\l) = \sum_{r=1}^{\l}
(-1)^{i_1+\cdots + i_{r-1}} c_1\otimes\cdots\otimes c_{r-1}
\otimes\delta_r(c_r)\otimes c_{r+1}\otimes\cdots\otimes c_\l,$$
where $\delta_r=\Delta^{(r)}$ if $i_r=1$ and
$\delta_r=\epsilon^{(r)}$ if $i_r=0$.

The simplest such instance is when all the exponents $d_i$
are equal to 1, in which case we have
$G=\Z^\l=\langle x_i\mid [x_i,x_j]=1 \rangle$.
Then $\widetilde{C}_\bullet(\Z^\l)$ is the usual free
$\Z\Z^\l$-resolution of $\Z$.  That is, $X_G$ is the
$\l$-torus $T^\l$, and $\widetilde{C}_\bullet=\widetilde{C}_\bullet(\Z^\l)$
is the equivariant augmented chain complex of the universal
cover of $T^\l$.   Specifically, $C_0=R=\Z\Z^\l$, $C_1=R^\l$
may be identified with
a free $R$-module with basis $\{e_1,\dots,e_{\l}\}$, and
$C_k = R^{\binom{\l}{k}} \cong \bigwedge^k C_1$.
With these identifications, the differential
$\Delta:C_k\to C_{k-1}$ may be expressed as
$$\Delta(e_{i_1}\wedge\dots\wedge e_{i_k}) = \sum_{r=1}^k
(-1)^{r-1}(x_{i_r}-1) (e_{i_1}\wedge\dots\wedge \hat
e_{i_r}\wedge\dots\wedge e_{i_k}).$$
\endexample

\subhead 3.4. IA-Products \endsubhead
An automorphism of a group $G$ is said to be an {\it IA-automorphism}
if it induces the identity automorphism on the abelianization
$H_1(G;\Z)=G/G'$.  The IA-automorphisms of $G$ form a (normal)
subgroup $\IA(G)$ of $\Aut(G)$.  The groups $\IA(F_n)$ have been much
studied; for example, it is known from the work of Nielsen and Magnus
that the natural map $\Aut(F_n)\to \GL(n,\Z)$ is surjective,
and that its kernel, $\IA(F_n)$, is finitely generated and
torsion-free. Let us say that a group $G$ is an {\it IA-product
of free groups} if it can be written as
$G=G_\l \rtimes_{\alpha_{\l}} \dots \rtimes_{\alpha_{3}}
G_2 \rtimes_{\alpha_{2}} G_1,$
with $G_q=F_{d_q}$, and $\alpha_q:G^{q-1}\to \IA(F_{d_q})$.

\proclaim{Proposition 3.5} Let $G$ be an IA-product of free groups.
Then the chain complex \linebreak
$C_\bullet(G)\otimes_G\Z$ has trivial boundary maps.
\endproclaim
\demo{Proof} It suffices to show that each of the matrices
$\Delta_{p_1,\dots,p_k}$ reduces to the zero matrix upon applying
the augmentation map $\epsilon:\Z G \to \Z$ to every entry.
This is accomplished by an inductive argument.\quad\qed
\enddemo
Let $b_j(G) = \rank H_j(G;\Z)$ be the $j^{\text{th}}$ Betti
number of $G$.  The following generalize results of Kohno \cite{27},
and Falk and Randell \cite{18}.
\proclaim{Corollary 3.6 (Factorization)} If $G$ is an IA-product
of free groups, then the homology groups of $G$ are torsion free,
and the Poincar\'e polynomial of $G$ factors into linear terms:
$$\sum_{j = 0}^\l b_j(G)t^j = \prod_{q=1}^\l (1 + d_qt).$$
\endproclaim
Let $\phi_k = \rank \Gamma_k(G)/\Gamma_{k+1}(G)$ denote the rank
of the $k^{\text{th}}$ lower central series quotient of $G$.
Noting that the group theoretic results of \cite{18} or \cite{30}
apply in our situation, we obtain:
\proclaim{Theorem 3.7 (LCS Formula)} If $G$ is an IA-product of
free groups, then in $\Z[[t]]$ we have
$$\sum_{j = 0}^\l b_j(G)(-t)^j = \prod_{q=1}^\l (1 - d_qt)
=\prod_{k\ge 1}(1-t^k)^{\phi_k}.$$
\endproclaim

\subhead{3.8. An Example}\endsubhead
We conclude this section with a detailed example of the
construction of the chain complex from section~2.   Recall that
$B_\l$ denotes the group of braids on $\l$ strings.  Let $B_\l^1$ be
the subgroup of braids that fix the endpoint of the last
string.  The group  $B_\l^1$ is the semidirect product of $B_{\l-1}$
with $F_{\l-1}$, determined by the Artin representation
$\alpha_{\l-1}:B_{\l-1}\to \Aut(F_{\l-1})$ (see~1.3).

\example{Example 3.9}  Let $G=B_4^1=F_3\rtimes_{\alpha_3} B_3$.  This
group admits the structure of an iterated semidirect product of free
groups, $G = F_3\rtimes_{\alpha_3} F_2\rtimes_{\mu_2} F_1$.
To see this, first note that
$B_3=\langle \sigma_1, \sigma_2 \mid \sigma_1\sigma_2\sigma_1=
\sigma_2\sigma_1\sigma_2 \rangle$ admits a semidirect product
structure $B_3=F_2\rtimes_{\mu_2} F_1$, coming from the Milnor
fibration of the discriminant singularity $\dd_3$ in $\C^3$ (see~7.5).
The group $F_1$ is generated by $x_{1,1}=\sigma_1$, the
group $F_2$ is generated by $x_{1,2}=\sigma_1\sigma_2^{-1},
x_{2,2}=\sigma_2^{-1}\sigma_1$, and the action
$\mu_2: F_1\to \Aut(F_2)$ is given by:
$$x_{1,1}:\cases x_{1,2}\mapsto x_{2,2}\\
x_{2,2}\mapsto x_{1,2}^{-1}x_{2,2}\endcases$$

Now let $F_3=\langle x_{1,3}, x_{2,3}, x_{3,3} \rangle$.
The Artin representation $\alpha_3: B_3 \to \Aut(F_3)$ is
given by:
$$x_{1,1}:
\cases
x_{1,3}\mapsto x_{1,3}x_{2,3}x_{1,3}^{-1}\\
x_{2,3}\mapsto x_{1,3}\\
x_{3,3}\mapsto x_{3,3}
\endcases
$$
$$
x_{1,2}:
\cases
x_{1,3}\mapsto x_{1,3}x_{3,3}x_{1,3}^{-1}\\
x_{2,3}\mapsto x_{1,3}\\
x_{3,3}\mapsto x_{3,3}^{-1}x_{2,3}x_{3,3}
\endcases
\qquad
x_{2,2}:
\cases
x_{1,3}\mapsto x_{1,3}x_{2,3}x_{1,3}^{-1}\\
x_{2,3}\mapsto x_{3,3}\\
x_{3,3}\mapsto x_{3,3}^{-1}x_{1,3}x_{3,3}
\endcases
$$

This establishes the iterated semidirect product structure on
$G=B_4^1$.  Carrying out the procedure described in section 2,
we see that the chain complex $C_\bullet(G)$ is given by
$$R^6 @>\pmatrix \Delta_{1,2,3}&\Delta_{2,3}&[\Delta_3]^2
\endpmatrix>> R^6\oplus R^3\oplus R^2
@>\pmatrix \Delta_{2,3}&[\Delta_3]^2&0\\
\Delta_{1,3}&0&\Delta_3\\ 0&\Delta_{1,2}&\Delta_2
\endpmatrix >>
R^3\oplus R^2\oplus R
@>\pmatrix \Delta_3\\ \Delta_2 \\ \Delta_1 \endpmatrix >> R$$
where $R=\Z G$, and
$$\Delta_1=\pmatrix x_{1,1}-1 \endpmatrix ,\qquad
\Delta_2=\pmatrix x_{1,2}-1\\x_{2,2}-1 \endpmatrix ,\qquad
\Delta_3=\pmatrix x_{1,3}-1\\x_{2,3}-1\\x_{3,3}-1 \endpmatrix, $$
$$\Delta_{1,2}=\pmatrix 1 & -x_{1,1}\\x_{1,1}x_{1,2}^{-1} &
1-x_{1,1}x_{1,2}^{-1} \endpmatrix ,\quad
\Delta_{1,3}=\pmatrix 1+(x_{1,3}-1)x_{1,1} & -x_{1,1}x_{1,3} & 0\\
-x_{1,1} & 1 & 0\\
0 & 0 & 1-x_{1,1} \endpmatrix ,$$
$$\Delta_{2,3}=\pmatrix
1+(x_{1,3}-1)x_{1,2} & 0 & -x_{1,2}x_{1,3}\\
-x_{1,2} & 1 & 0\\
0 & -x_{1,2}x_{3,3}^{-1} & 1-x_{1,2}x_{3,3}^{-1}(x_{2,3}-1) \\
1+(x_{1,3}-1)x_{2,2} & -x_{2,2}x_{1,3} & 0\\
0 & 1 & -x_{2,2}\\
-x_{2,2}x_{3,3}^{-1} & 0 & 1-x_{2,2}x_{3,3}^{-1}(x_{1,3}-1)
\endpmatrix ,$$
and $\DS{\Delta_{1,2,3}=\pmatrix -I_3 & I_3 - \Delta_{1,3}\\
-A & -I_3 + A \endpmatrix}$, where
$$A=\pmatrix (1-x_{1,3})x_{1,1}x_{1,2}^{-1} &
(1-x_{1,3})x_{1,1}x_{1,2}^{-1}x_{1,3} &
x_{1,1}x_{1,2}^{-1}x_{1,3}x_{2,3}\\
0 & x_{1,1}x_{1,2}^{-1} & 0 \\
x_{1,1}x_{1,2}^{-1}x_{2,3}^{-1} &
x_{1,1}x_{1,2}^{-1}x_{2,3}^{-1}(x_{1,3}-1) & 0
\endpmatrix .
$$
\remark{Remark 3.10}
This chain complex can be used to compute the homology of $G$
with various (twisted) coefficients. For example, the homology
$H_*=H_*(G;\Z)$ with trivial $\Z$ coefficients is given by
$H_0=\Z, H_1=\Z^2, H_2=\Z^2, H_3=\Z$.  Note that the monodromy
automorphisms $\alpha_2$ and $\alpha_3$ defining $G$ are {\it not}
IA-automorphisms, and that the Factorization formula from
Corollary 3.6 does not hold for this group.
\endremark
\endexample

\head 4. Generalized Magnus Representations
\endhead

Let $G$ be a group that admits the structure of an iterated
semidirect product of finitely generated free groups.
In this section, we use the construction of the chain complex
$C_\bullet(G)$ to define representations of groups $\Gamma$
which act compatibly on $G$.

\definition{Definition 4.1}  An automorphism $\psi\in \Aut(G)$
is said to be {\it compatible} with the iterated semidirect product
structure $G=G_\l\rtimes_{\alpha_\l}\cdots \rtimes_{\alpha_3}
G_2 \rtimes_{\alpha_2} G_1$ if it satisfies the following conditions:
\roster
\item $\forall g \in G_p$, we have $\psi(g) \in G_p$; and
\smallskip
\item $\forall g_p \in G_p$, $g_q \in G_q$, $p<q$, we have
$\alpha_q^p(\psi(g_p))(\psi(g_q)) =
\psi(\alpha_q^p(g_p)(g_q))$ in $G_q$.
\endroster
\enddefinition

Given a compatible automorphism $\psi$ of $G$, let
$\psi_p\in\Aut(G_p)$ be the restriction of $\psi$ to $G_p$
(well defined by condition (1) above).  Condition (2) can be
restated as:
\roster
\item"($2'$)" $\forall x\in G_p$, $p<q$, we have
$\alpha_q^p(\psi_p(x)) =
\psi_q\circ \alpha_q^p(x)\circ \psi_q^{-1}$ in $\Aut(G_q)$.
\endroster

It is easily checked that the compatible automorphisms of $G$
form a subgroup of $\Aut(G)$, which we shall denote by
$\Aut^{\rtimes}(G)$.  The compatibility conditions are quite
restrictive.  For example, if $\DS{G=\times_1^\l G_i}$, then
$\DS{\Aut^{\times}(G)=\times_1^\l \Aut(G_i)}$; in particular,
$\Aut^{\times}(\Z^\l) \cong (\Z_2)^\l$.

\subhead{4.2}\endsubhead  Now assume $G_p=F_{d_p}$, for each $p$,
$1\le p\le\l$, and let $C_\bullet=C_\bullet(G)$ be the chain complex
of free modules over $R=\Z G$ constructed in section 2.  Any
automorphism $\psi:G\to G$ gives rise to a chain equivalence
$\Psi_\bullet:C_\bullet\to C_\bullet$.  We shall
explicitly describe the chain map $\Psi_\bullet$ in the case where
$\psi$ belongs to $\Aut^{\rtimes}(G)$.

Let $J(\psi_p)\in\Aut_{\Z G_p}((\Z G_p)^{d_p})$ be the Jacobian matrix
of Fox derivatives of the automorphism $\psi_p\in\Aut(G_p)$.
Let $J_p(\psi):=\id_R\otimes J(\psi_p)\in\Aut_R(R^{d_p})$
be the automorphism obtained from  $J(\psi_p)$ by extension
of scalars.  Define the {\it higher-order Jacobians} of $\psi$
recursively as follows:
$$J_{p_1,\dots,p_k}(\psi):=
\left[J_{p_k}(\psi)\right]^{d_{p_1}\cdots d_{p_{k-1}}}\cdot
\tilde\rho_{p_k}(J_{p_1,\dots,p_{k-1}}(\psi)):
R^{d_{p_1}\cdots d_{p_k}}\to R^{d_{p_1}\cdots d_{p_k}}.$$

Now define the map $\Psi_k: C_k\to C_k$ by specifying its restriction
to the summand $R^{d_{p_1}\cdots d_{p_k}}$ of $C_k$ to be the composition
of the higher-order Jacobian $J_{p_1,\dots,p_k}(\psi)$ with the extension
$\tilde\psi:R^{d_{p_1}\cdots d_{p_k}}\to R^{d_{p_1}\cdots d_{p_k}}$:
$$\Psi_k=\bigoplus_{1\le p_1<\cdots < p_k\le\l}
J_{p_1,\dots,p_k}(\psi)\circ\tilde\psi.$$
(The map $\Psi_0: C_0\to C_0$ is just $\tilde\psi:R\to R$.)
In order to prove that $\Psi_{\bullet}: C_{\bullet}\to C_{\bullet}$
is a chain map, we first need a lemma.

\proclaim{Lemma 4.3} For all $x\in G_p$ and $p<q$, we have
$\rho_q(\psi_p(x))=J_q(\psi)^{-1} \cdot \tilde\psi(\rho_q(x))
\cdot J_q(\psi)$.
\endproclaim
\demo{Proof}  Compute:
$$\allowdisplaybreaks
\alignat2
\rho_q(\psi_p(x))
&=\psi_p(x)\cdot J(\alpha_q^p(\psi_p(x)))
&&\qquad\text{by 2.7}\\
&=\psi_p(x)\cdot J(\psi_q\circ\alpha_q^p(x)\circ\psi_q^{-1})
&&\qquad\text{by ($2'$) of 4.1}\\
&=\psi_p(x)\cdot\tilde\psi_q\circ\widetilde{\alpha_q^p(x)}
(J(\psi_q^{-1}))\cdot J(\psi_q\circ\alpha_q(x))
&&\qquad\text{by 2.3}\\
&=\psi_p(x)\cdot\tilde\psi_q\circ\widetilde{\alpha_q^p(x)}\circ
\tilde\psi_q^{-1}(J(\psi_q)^{-1})\cdot
J(\psi_q\circ\alpha_q^p(x))
&&\qquad\text{by 2.3}\\
&=\psi_p(x)\cdot\widetilde{\alpha_q^p(\psi_p(x))}(J_q(\psi)^{-1})
\cdot J(\psi_q\circ\alpha_q^p(x))
&&\qquad\text{by ($2'$) of 4.1}\\
&=\psi_p(x)\cdot\psi_p(x)^{-1}\cdot J_q(\psi)^{-1}\cdot\psi_p(x)
\cdot J(\psi_q\circ\alpha_q^p(x))
&&\qquad\text{by \thetag{1.2}}\\
&=J_q(\psi)^{-1}\cdot\psi(x)
\cdot \tilde\psi(J(\alpha_q^p(x)))\cdot J_q(\psi)
&&\qquad\text{by 2.3}\\
&=J_q(\psi)^{-1}\cdot \tilde\psi(\rho_q(x))\cdot J_q(\psi)
&&\qquad\text{by 2.7}\qquad\quad\qed
\endalignat
$$
\enddemo

\proclaim{Proposition 4.4} The map $\Psi_\bullet:C_\bullet\to
C_\bullet$ is a chain equivalence.
\endproclaim
\demo{Proof}  To show that $\Psi_\bullet:C_\bullet\to C_\bullet$
is a chain map, we must show that
$\Psi_{k-1} \circ \Delta = \Delta \circ \Psi_k$, for $1\le k\le \l$.
We accomplish this by induction on $k$, with
the case $k=1$ following from diagram \thetag{2.2}.

By virtue of the direct sum decompositions of $C_\bullet$ and
$\Psi_\bullet$, it is enough to show that diagrams of the form
$$\CD
R^{d_{p_1}\cdots d_{p_k}} @>\Delta^{p_1,\dots,p_k} >>
\bigoplus_{i=1}^k R^{d_{p_1}\cdots \hat d_{p_i}\cdots
d_{p_k}}\\
@VV\Psi_k V     @VV\Psi_{k-1} V\\
R^{d_{p_1}\cdots d_{p_k}} @>\Delta^{p_1,\dots,p_k} >>
\bigoplus_{i=1}^k R^{d_{p_1}\cdots \hat d_{p_i}\cdots
d_{p_k}}\\
\endCD$$
commute.  Since
$$\Delta^{p_1,\dots,p_k} =
\left(\Delta_{p_1\dots,p_{k}},\left[\Delta_{p_2,\dots,p_k}\right]^{d_{p_1}},
\dots,\left[\Delta_{p_i,\dots,p_k}\right]^{d_{p_1}\cdots d_{p_{i-1}}},
\dots,\left[\Delta_{p_k}\right]^{d_{p_1}\cdots d_{p_{k-1}}}\right),$$
this amounts to showing that
$$J_{p_1,\dots,p_{k}}(\psi)\cdot
\left[\Delta_{p_i,\dots,p_k}\right]^{d_{p_1}\cdots d_{p_{i-1}}} =
\tilde\psi\left(
\left[\Delta_{p_i,\dots,p_k}\right]^{d_{p_1}\cdots d_{p_{i-1}}}\right)
\cdot J_{p_1,\dots,\hat p_i,\dots,p_k}(\psi)$$
for $1\le i \le k$. These equalities all follow from the definitions using
induction, together with Lemma 4.3.

To complete the proof, we must show that the
chain map $\Psi_\bullet$ is, in fact, a chain equivalence.
This follows directly from the definitions.\quad\qed
\enddemo

\example{Example 4.5}  The simplest example of this construction
is in the case of a direct product $G=\times_1^{\l} G_p$, where
$G_p=F_{d_p}$.  Using the decomposition $C_k(G) =\bigoplus
C_{i_1}(G_1)\otimes\cdots \otimes  C_{i_\l}(G_\l)$ from 3.2,
we can write the chain map induced by $\psi=\times_1^\l \psi_p$ as
$\Psi_k = \bigoplus \Psi_{i_1}^{(1)}\otimes \cdots \otimes
\Psi_{i_\l}^{(\l)}$,
where $\Psi_{i_r}^{(r)}=J(\psi_r)\circ\tilde\psi_r$ if $i_r=1$ and
$\Psi_{i_r}^{(r)}=\tilde\psi_r$ if $i_r=0$.
\endexample

\example{Example 4.6}   We further illustrate the construction
using the notations and computations of Example 3.9.
Recall the group $G=B_4^1=F_3\rtimes_{\alpha_3}
F_2\rtimes_{\mu_2} F_1$.  Consider the normal subgroup
$\Cal G=F_3\rtimes_{\alpha_3} F_2$,
where $F_2=\langle x_{1,2}, x_{2,2}\rangle$,
$F_3=\langle x_{1,3}, x_{2,3}, x_{2,3}\rangle$,
and define the automorphism $\psi\in\Aut^{\rtimes}(\Cal G)$
to be conjugation by $x_{1,1}\in F_1$.

The chain complex $C_\bullet(\Cal G)$ is of the form
$$\Cal R^6 @>{\pmatrix \Delta_{2,3}&[\Delta_3]^2\endpmatrix}>>
\Cal R^3\oplus \Cal R^2
@>{\pmatrix \Delta_3 \\ \Delta_2 \endpmatrix}>> \Cal R$$
where $\Cal R=\Z \Cal G$, and the boundary maps are given by
restricting the boundary maps of $C_\bullet(G)$.
Carrying out the construction described in 4.2, the chain
equivalence $\Psi_\bullet:C_\bullet(\Cal G)\to C_\bullet(\Cal G)$
induced by $\psi$ can be written as:
$$\CD
\Cal R^6 @>>> \Cal R^3\oplus \Cal R^2 @>>> \Cal R\\
@VV{J_{2,3}\circ\tilde\psi}V @VV{(J_3\oplus J_2)\circ\tilde\psi}V
@VV{\tilde\psi}V\\
\Cal R^6 @>>> \Cal R^3\oplus \Cal R^2 @>>> \Cal R\\
\endCD$$
where $J_2=J(\left.\psi\right|_{F_2})$, $J_3=J(\left.\psi\right|_{F_3})$,
and $J_{2,3}=\left[ J_3\right]^2\cdot\tilde\rho_2(J_2)$ are given by
$$
J_2 = x_{1,1}^{-1}\cdot (I_2 - \Delta_{1,2})
=\pmatrix  0 &1 \\ -x_{1,2}^{-1} & x_{1,2}^{-1}\endpmatrix ,
$$
$$
J_3 = x_{1,1}^{-1}\cdot (I_3 - \Delta_{1,3}) ,
\qquad J_{2,3} = x_{1,1}^{-1}\cdot (I_6 - \Delta_{1,2,3}) .
$$
\endexample

\subhead{4.7}\endsubhead  In order to proceed with the construction
of generalized Magnus representations, we need to see how the
higher-order Jacobians behave under composition.  The following
result can be viewed as a generalization of the chain rule of
Fox Calculus (Lemma 2.3).  Let $G=\rtimes_{p=1}^{\l} F_{d_p}$,
and consider $\phi, \psi\in\Aut^{\rtimes}(G)$.

\proclaim{Proposition 4.8 (Chain Rule)}
$J_{p_1,\dots,p_k}(\psi\circ \phi)=
\tilde\psi(J_{p_1,\dots,p_k}(\phi))\cdot J_{p_1,\dots,p_k}(\psi)$.
\endproclaim

\demo{Proof}  This is proved by induction on $k$, with the case
$k=1$ following from Lemma 2.3.

For the inductive step we need Lemma 4.3, which, we recall,
states that for $x\in G^{q-1}$, we have
$\rho_q(\psi(x)) = J_q(\psi)^{-1}\cdot \tilde\psi(\rho_q(x))
\cdot J_q(\psi)$.  It follows immediately that
$\tilde\rho_q\left(\tilde\psi\left(\sum n_x x\right)\right) =
J_q(\psi)^{-1}\cdot \tilde\psi\left(\tilde\rho_q(\sum n_x x)\right)
\cdot J_q(\psi)$.
If $A$ is a $d \times d$ matrix with entries in $R$,
this implies that
$$\tilde\rho_q\left(\tilde\psi(A)\right) =
\left[ J_q(\psi)^{-1}\right]^d\cdot
\tilde\psi(\tilde\rho_q(A))\cdot \left[J_q(\psi)\right]^d.$$

Using the definition of higher-order Jacobians, the case $k=1$,
the inductive hypothesis, the fact that $\tilde\rho_q$ is a
homomorphism (see 2.8), and the above formula, we get:
$$
\align
J_{p_1,\dots,p_k}(\psi\circ \phi)
&=\left[J_{p_k}(\psi\circ\phi)\right]^{d_{p_1}\cdots d_{p_{k-1}}}\cdot
\tilde\rho_{p_k}(J_{p_1,\dots,p_{k-1}}(\psi\circ\phi))\\
&=\left[\tilde\psi(J_{p_k}(\phi))\cdot J_{p_k}(\psi)
\right]^{d_{p_1}\cdots d_{p_{k-1}}}
\cdot \tilde\rho_{p_k}\left( \tilde\psi(J_{p_1,\dots,p_{k-1}}(\phi))
\cdot J_{p_1,\dots,p_{k-1}}(\psi) \right)\\
&=\left[\tilde\psi\left(J_{p_k}(\phi)\right)\right]^{d_{p_1}\cdots d_{p_{k-1}}}
\cdot\left[J_{p_k}(\psi)\right]^{d_{p_1}\cdots d_{p_{k-1}}}
\cdot\tilde\rho_{p_k}\left(\tilde\psi(J_{p_1,\dots,p_{k-1}}(\phi))\right)\\
&\qquad\qquad\cdot\tilde\rho_{p_k}\left(J_{p_1,\dots,p_{k-1}}(\psi)\right)\\
&=\tilde\psi\left(\left[J_{p_k}(\phi)\right]^{d_{p_1}\cdots d_{p_{k-1}}}\right)
\cdot\tilde\psi\left(\tilde\rho_{p_k}(J_{p_1,\dots,p_{k-1}}(\phi))\right)
\cdot\left[J_{p_k}(\psi)\right]^{d_{p_1}\cdots d_{p_{k-1}}}\\
&\qquad\qquad\cdot\tilde\rho_{p_k}\left(J_{p_1,\dots,p_{k-1}}(\psi)\right)\\
&=\tilde\psi(J_{p_1,\dots,p_k}(\phi))\cdot J_{p_1,\dots,p_k}(\psi)
\qquad\quad\qed
\endalign
$$
\enddemo

\subhead{4.9}\endsubhead   Now consider a group $\Gamma$
that acts compatibly on $G$.  That is, we are given
a homomorphism $\Phi:\Gamma\to \Aut^{\rtimes}(G)$.  For each
$\gamma\in\Gamma$, the construction of 4.2 yields a chain map
$\Phi(\gamma)_{\bullet}: C_\bullet(G) \to C_\bullet(G)$.  But
these chain maps are not $R$-linear in general, so a judicious
extension of scalars is required in order to define a
representation of $\Gamma$.  The idea is suggested by the
original approach followed by Magnus \cite{36} (see \cite{5}, p.~115).

\definition{Definition 4.10} A homomorphism $\tau: G\to K$
is said to be {\it $\Phi$-invariant} if
$\tau(\Phi(\gamma)(g)) = \tau(g)$ for all
$g\in G$ and $\gamma\in\Gamma$.
\enddefinition

Let $R=\Z G$, $S=\Z K$, $\tilde\tau:R\to S$ the extension
of $\tau$ to group rings, and $S\otimes_R -$ the
extension of scalars functor defined by $\tilde{\tau}$.
Applying this functor, we obtain a chain complex of free $S$-modules,
$S\otimes_R C_{\bullet}(G)$.  We are now ready to state
the main theorem of this section.

\proclaim{Theorem 4.11} Suppose $G=\rtimes_{p=1}^{\l} G_p$ is an
iterated semidirect product of free groups,
$\Phi:\Gamma\to \Aut^{\rtimes}(G)$ is a compatible action of a group
$\Gamma$ on $G$,  and $\tau:G\to K$ is a $\Phi$-invariant
homomorphism.  Given $\gamma\in \Gamma$, let
$\Phi^{\tau}_{\bullet}(\gamma)=\id_S\otimes \Phi(\gamma)_{\bullet}:
S\otimes_R C_{\bullet}(G) \to S\otimes_R C_{\bullet}(G)$.
Then, for each $k$, $1\le k\le \l$,
\roster
\item"(i)" the map $\Phi^{\tau}_k(\gamma):S\otimes_R C_k(G)\to
S\otimes_R C_k(G)$ is S-linear;
\item"(ii)"  the map $\Phi^{\tau}_k(\gamma)$ is a chain equivalence; and
\item"(iii)"  the map $\Phi^{\tau}_k:\Gamma\to\Aut_S(S\otimes_R C_k(G)),
\gamma \mapsto \Phi_k^{\tau}(\gamma)$, is a homomorphism.
\endroster
\endproclaim
\demo{Proof}  (i) For a fixed $\gamma\in\Gamma$, write
$\psi=\Phi(\gamma)\in \Aut^{\rtimes}(G)$, and let
$\Psi_\bullet:C_\bullet\to C_\bullet$ be the chain map
induced by $\psi$.  Recall that $\Psi_0:C_0 \to C_0$ is
identified with $\tilde\psi: R\to R$.  Recall also that
$\Psi_k:C_k \to C_k$ is the composition of a certain $R$-linear
map with the non-linear map $\tilde\psi: C_k \to C_k$, which
is a direct sum of copies of $\tilde\psi: R\to R$.
Thus, it is enough to show that
$\Phi^{\tau}_0=\id_S\otimes\tilde\psi: S\otimes_R R \to S\otimes_R R$
is an $S$-linear map.  We will show more, namely
$$\id_S\otimes \tilde\psi = \id_{S\otimes_R R}.\tag{$\ddagger$}$$

Let $\omega:S\otimes_R R @>\sim>> S$ be the canonical isomorphism
given by  $\omega(s\otimes r) = s\tilde\tau(r)$.
The $\Phi$-invariance condition, $\tau\circ\psi=\tau$, yields:
$$\omega(\Phi^{\tau}_0(s\otimes r)) = \omega(s\otimes \tilde\psi(r)) =
s\tilde\tau(\tilde\psi(r)) = s\tilde\tau(r)=\omega(s\otimes r),$$
proving the claim.

\roster
\item"(ii)" This follows from Proposition 4.4.
\item"(iii)" This follows from Proposition 4.8 and \thetag{$\ddagger$}.
\quad\qed
\endroster
\enddemo

In the special case where $G=F_n$ and $\Gamma < \Aut(F_n)$,
representations such as the above were introduced by Magnus \cite{36}
(see \cite{5} for details).  We therefore refer to the
homomorphisms $\Phi^{\tau}_k:\Gamma\to\Aut_S(S\otimes_R C_k(G))$ as
{\it generalized Magnus representations}.
Since the maps $\Phi^{\tau}_\bullet(\gamma)$ are chain maps,
we also obtain {\it homological Magnus representations}
$\bar\Phi^{\tau}_k:\Gamma \to \Aut_S H_k(S\otimes_R C_{\bullet}(G))$.
However, note that these homology groups need not be free $S$-modules
in general.  In such a situation, one may still be able to ``reduce''
$\Phi^{\tau}_k$ by restricting to a free, invariant submodule of
$S\otimes_R C_k(G)$---see e.g.~5.9.

\remark{Remark 4.12} There is an alternate way to interpret
these representations.  Recall that, in forming the tensor
product $S\otimes_R C_{\bullet}(G)$, we view
$S$ as a right $R$-module via $s\cdot r=s\tilde\tau(r)$.
Using the involution of the group ring $R=\Z G$ given by
$\overline{\sum n_g g} = \sum n_g g^{-1}$, we can turn
$S$ into a left $R$-module by setting
$r\cdot s=s\tilde\tau(\bar r)$, and
form the tensor product $C_{\bullet}(G) \otimes_G S$.
We then have a chain equivalence
$S\otimes_R C_{\bullet}(G) @>\sim>>C_{\bullet}(G) \otimes_G S$
given by $s\otimes c\mapsto c\otimes s$, inducing an
isomorphism
between $H_\bullet(S\otimes_R C_{\bullet}(G))$ and
$H_{\bullet}(C_{\bullet}(G)\otimes_G S)$.
Thus, we can view the generalized Magnus representations as
$\Phi^{\tau}_k:\Gamma \to \Aut_S(C_k(G)\otimes_G S)$, respectively
$\bar\Phi^{\tau}_k:\Gamma \to \Aut_S H_k(G; S)$.  When $K$ is an
abelian group, the coefficients module $S=\Z K$ is determined
by the representation $\hat\tau: G\to \Aut S$, given by
$\hat\tau(g)(s) = \tau(g^{-1})s$.
\endremark

\example{Example 4.13}  The simplest example is where
$\tau: G\to \{1\}$ is the trivial homomorphism.  If
$\Gamma$ acts compatibly on $G$, the resulting representations,
$\Phi^{\tau}_k:\Gamma \to \Aut(C_k(G)\otimes_G \Z)$ and
$\bar\Phi^{\tau}_k:\Gamma \to \Aut H_k(G; \Z)$, can be non-trivial,
even when $G=F_n$, see \cite{5}, p.~117.  However, if $\Gamma$ acts
by $\IA$-automorphisms of $G$, then obviously $\bar\Phi^{\tau}_k$
is trivial.
\endexample

\example{Example 4.14}  Let $\Gamma$ be a group that acts compatibly
on $G$, and assume that the action $\Phi$ of $\Gamma$ on $G$
factors as
$\Phi:\Gamma\to \Aut^{\rtimes}(G)\cap\IA(G) \hookrightarrow
\Aut^{\rtimes}(G)$.
Then the abelianization map $\ab:G \twoheadrightarrow G/G'$, and more
generally, maps of the form $\tau:G @>{\ab}>> G/G'
\twoheadrightarrow K$ (where $K$ is abelian) are $\Phi$-invariant.
See 5.8 and 5.9 for examples of generalized Magnus
representations obtained this way.
\endexample

\head 5. Representations of Braid Groups
\endhead

We use the techniques developed above to define new
linear representations of braid groups.   Detailed discussion
of these representations is deferred to \cite{12}.

\subhead{5.1. Generalized Burau Representations}
\endsubhead
For $1\le \l< n$, let $P_{n,\l}=\ker (P_n\to P_\l)$
denote the kernel of the homomorphism from $P_n$ to $P_\l$
defined by forgetting the last $n-\l$ strands.  Then
$P_{n,\l}=\rtimes_{p=\l}^{n-1} F_p$ is generated by
$\{A_{i,j}\}$ with $\l< j\le n$ and $1\le i<j$.  (Note that $P_n=P_{n,1}$.)
The braid group $B_\l$ acts on $P_{n,\l}$ in a natural fashion.
On each free factor $F_p=F_{\l}*F_{p-\l}$ of $P_{n,\l}$, it acts
by the Artin representation on $F_{\l}$ (see~1.3), and acts trivially
on $F_{p-\l}$. It is readily checked that the action so defined,
$\Phi_{\l,n}: B_\l\to\Aut(P_{n,\l})$, is compatible
with the iterated semidirect product structure of $P_{n,\l}$.
The semidirect product $P_{n,\l}\rtimes_{\Phi_{\l,n}} B_\l$
is the group $B_n^{n-\l}$ of braids that fix the endpoints of
the last $n-\l$ strings.

Let $\Z=\langle t \rangle$, and identify the group
ring, $\Z\Z$, with the ring of Laurent polynomials in $t$,
$\Lambda=\Z[t, t^{-1}]$.   Fix $m$, $1\le m\le n-\l$, and define a
homomorphism $\tau:P_{n,\l}\to\Z$ by
$\tau(A_{r,s}) =  t$ if $n-m+1\le s\le n$, and $\tau(A_{r,s})=1$
if $\ell+1\le s\le n-m$.
Since $\tau(A_{r,s})=\tau(A_{p,q})$ if $s=q$, the homomorphism $\tau$ is
invariant with respect to the action $\Phi_{\l,n}$ of $B_\l$ on $P_{n,\l}$.
We thus obtain by Theorem~4.11
generalized Magnus
representations of the braid group,
$$\align
\beta^m_{\l,n-\l,k}=(\Phi_{\l,n})^{\tau}_k&: B_\l\to\Aut_{\Lambda}
\left(\Lambda \otimes_{\Z P_{n,\l}} C_k(P_{n,\l})\right),\\
\bar\beta^m_{\l,n-\l,k}&:B_\l\to\Aut_{\Lambda}
H_k(\Lambda \otimes_{\Z P_{n,\l}} C_\bullet(P_{n,\l}))
\endalign$$
for each $k$, $1\le k\le n-\l$.
If $\l=n-1$, we have $P_{n,\l}=F_\l$, $k=1$, and
the action $\Phi_{\l,\l+1}=\alpha_{\l}$ is the Artin representation.
In this instance, it is easy to see that
$\beta_{\l}=\beta^1_{\l,1,1}:B_\l\to \GL(\l,\Lambda)$
is the Burau representation, and
$\bar\beta_{\l}:B_\l\to \GL(\l-1,\Lambda)$
is the reduced Burau representation (see~\cite{5}, \cite{26},
\cite{34}).

\example{Example 5.2}
In the case $\l=3$, $n=5$, $m=1$, we have
$H_2(\Lambda \otimes_{\Z P_{5,3}} C_\bullet(P_{5,3}))=\Lambda^6$.
Thus the above construction yields
a generalized Burau representation
$\bar\beta_{3,2,2}^1:B_3\to\GL(6,\Lambda)$.
This representation is given by
$$\allowdisplaybreaks
\align
\sigma_1\mapsto
&\pmatrix
0&0&-t-t^2&t&0&0\\
0&0&-1&0&0&0\\
1-t&-1&0&0&0&0\\
1&-1-t&0&0&0&0\\
0&0&0&0&-t&0\\
0&0&0&0&-t-t^2&1
\endpmatrix,\\
\\
\sigma_2\mapsto
&\pmatrix
1&-t-t^2&0&0&0&0\\
0&-t&0&0&0&0\\
0&0&0&0&-t&0\\
0&0&0&0&-t^2-t^3&t^2\\
0&0&-t^{-1}&t^{-1}-1&0&0\\
0&0&-1-t^{-1}&t^{-1}&0&0
\endpmatrix.
\endalign$$
The characteristic polynomial of each of these matrices is
$(\lambda-1)^2(\lambda+1)(\lambda+t)(\lambda^2-t^3)$.
It follows that, unlike the Burau representation,
this representation does not factor through the Hecke algebra
$H(3,t)$, see \cite{26}.
\endexample

\remark{Remark 5.3} In \cite{29}, Kohno uses a vanishing theorem
involving the groups $P_{n,\l}$ to construct representations of
the braid group generalizing the (reduced) Burau representation.
See section~6 for a generalization of this vanishing theorem.
Like the above representations, those generated by Kohno do not,
in general, factor through the Hecke algebra.  The construction
of the generalized Magnus representations presented here was,
in part, motivated by Kohno's work.
\endremark

\remark{Remark 5.4}  Our generalized Burau representations are
powerful enough to detect braids which belong to the kernel of the
Burau representation.  Answering a long-standing question,
Moody \cite{40} showed that $\beta_{\l}:B_\l\to\GL(\l,\Lambda)$
is not faithful for $\l\ge 9$.  This was sharpened to $\l\ge 6$
by Long and Paton \cite{34}, who found the following braid in
$\ker(\beta_{6})$:  $\xi=[\zeta^{-1}\sigma_5 \zeta, (\sigma_2
\sigma_3\sigma_4\sigma_5)^5]$, where $\zeta=\sigma_4^{-1}\sigma_5
\sigma_3^{-1}\sigma_4\sigma_2^{-1}\sigma_3^{-3}\sigma_1^3\sigma_5
\sigma_4\sigma_3^{-1}\sigma_2^{-1}\sigma_1$.  A computation
shows that the representation $\bar\beta^1_{6,2,2}:B_6\to\GL(30,\Lambda)$
detects the braid $\xi$: $\bar\beta^1_{6,2,2}(\xi) \ne I_{30}$.
\endremark

\subhead{5.5}\endsubhead
The above construction may be applied in a more general
context so as to yield representations of the braid group
with several parameters.  Fix $m$, $1\le m\le n-\l$, and consider
the free abelian group $\Z^m$ generated by
$\{t_s\mid n-m+1\le s\le n\}$.  Identify the group ring $\Z\Z^m$
with the ring $\Lambda_m$ of Laurent polynomials in the variables
$\{t_s\}$.  Let $\tau:P_{n,\l}\to\Z^m$ be the homomorphism
defined by $\tau(A_{r,s}) = t_s$ if $n-m+1\le s\le n$ and
$\tau(A_{r,s})=1$ otherwise.  As before, the
homomorphism $\tau$ is invariant with respect to the action
$\Phi_{\l,n}: B_\l\to\Aut^{\rtimes}(P_{n,\l})$.
Thus we obtain generalized Magnus representations,
$$\align
\eta^m_{\l,n-\l,k}=(\Phi_{\l,n})^{\tau}_k&: B_\l\to\Aut_{\Lambda_m}
\left(\Lambda_m\otimes_{\Z P_{n,\l}} C_k(P_{n,\l})\right),\\
\bar\eta^m_{\l,n-\l,k}&: B_\l\to\Aut_{\Lambda_m}
H_k(\Lambda_m\otimes_{\Z P_{n,\l}} C_\bullet(P_{n,\l}))
\endalign$$
for each $k$,
$1\le k\le n-\l$, depending on $m$ parameters. Clearly,
$\eta^1_{\l,q,r}=\beta^1_{\l,q,r}$ and
$\bar\eta^1_{\l,q,r}=\bar\beta^1_{\l,q,r}$.

\example{Example 5.6} In the case $\l=3$, $n=5$, $m=2$, we have
$H_2(\Lambda_2\otimes_{\Z P_{5,3}}
C_\bullet(P_{5,3}))=\Lambda_2^6$. We therefore obtain a representation
$\bar\eta_{3,2,2}^2: B_3\to\GL(6,\Lambda_2)$.
This representation is given below, where we denote the generators of
$\Z^2$ by
$s$ and $t$ (as opposed to $t_4$ and $t_5$) to simplify notation.
$$\allowdisplaybreaks
\align
\sigma_1 \mapsto&
\pmatrix
0&0&0&s&-s&-s\\
0&1&0&1 + st&0&-1\\
0&0&0&-t&0&0\\
-st&0&-s^2&-t&-s^2t&s^2 - s^2t\\
0&0&1&-1&1&t\\
0&0&s+ st&-t - t^2&0&-s
\endpmatrix\\
\\
\sigma_2 \mapsto&
\pmatrix
-t - st&st&-s^2 + s^2t&0&0&0\\
-t&0&-s + st&0&0&0\\
t^2&0&-s + st&0&st&0\\
s&0&-1&1&s^2&0\\
t&0&s&0&0&0\\
-1 - t&0&-1 - s - t&0&-s&1
\endpmatrix
\endalign
$$
If one sets $s$ equal to $t$ in the above representation
(i.e.~applies the map $\tilde\phi:\Lambda_2\to\Lambda$ defined by
$\phi(s)=t$, $\phi(t)=t$), the resulting representation
$B_3\to\GL(6,\Lambda)$ is precisely the one parameter representation
$\bar\beta_{3,2,2}^2:B_3\to\GL(6,\Lambda)$ defined in 5.1.
\endexample
\remark{Remark 5.7}
Methods for ``lifting'' representations from $B_{\l+1}$ to $B_{\l}$ are
given by L\"udde and Toppan in \cite{35}, and by Birman, Long, and
Moody in \cite{6}.  Either of these techniques may be used to generate
braid group representations which depend on several parameters.  For
instance, L\"udde and Toppan obtain an $m$ parameter representation
$\nu_{\l}^m$ of $B_{\l}$ by successively lifting the trivial
representation of $B_{\l+m}$.  We have checked that the (reduced)
representation $\bar\nu_3^2$ is equivalent to $\bar\eta_{3,2,2}^2$.
We conjecture that $\nu_{\l}^m\cong\eta_{\l,m,m}^m$ and
$\bar\nu_{\l}^m\cong\bar\eta_{\l,m,m}^m$ for all $\l$ and $m$.
\endremark

\subhead{5.8. Generalized Gassner Representations}
\endsubhead
For each $n>\l$, the pure braid group $P_\l$ acts on the group
$P_{n,\l}=F_{n-1}\rtimes\cdots\rtimes F_{\l}$ by restriction
of the action $\Phi_{\l,n}$.  This (compatible) action is the
``usual'' one, discussed  in Example~1.3 (see also section~6).
The semidirect product $P_{n,\l}\rtimes_{\Phi_{\l,n}} P_\l$ is,
of course, the (entire) pure braid group $P_n$.

Fix $m$, $1\le m\le n-\l$, and let $N=\binom{n}2-\binom{n-m}2$.
Consider the free abelian group $\Z^N$ generated by
$\{t_{r,s}\mid 1\le r<s, n-m+1\le s\le n\}$, and let
$\tau:P_{n,\l}\to\Z^N$ be the homomorphism defined by
$\tau(A_{r,s}) = t_{r,s}$ if $n-m+1\le s\le n$
and $\tau(A_{r,s}) = 1$ otherwise.  Checking that $\tau$
is invariant with respect to the action $\Phi_{\l,n}$
of $P_\l$ on $P_{n,\l}$, we obtain generalized Magnus
representations of the pure braid group,
$$\align
\theta^m_{\l,n-\l,k}=(\Phi_{\l,n})^{\tau}_k&:
P_\l\to\Aut_{\Lambda_N}
\left(\Lambda_N\otimes_{\Z P_{n,\l}} C_k(P_{n,\l})\right),\\
\bar\theta^m_{\l,n-\l,k}&:P_\l\to\Aut_{\Lambda_N}
H_k(\Lambda_N\otimes_{\Z P_{n,\l}} C_\bullet(P_{n,\l}))
\endalign$$
for each $k$,
$1\le k\le n-\l$. In the special case $\l=n-1$ (and thus $m=1$), we have
$N=\l$ and
$P_{n,\l}=F_\l$.  In this instance, it is easy to see that
$\theta_{\l}=\theta_{\l,1,1}^1:P_\l\to \GL(\l, \Lambda_{\l})$
is the Gassner representation.  Note however that
$\bar\theta_{\l}:P_\l\to\Aut_{\Lambda_\l} H_1(
\Lambda_\l\otimes_{\Z F_{\l}} C_\bullet(F_{\l}))$ is not the reduced
Gassner representation for $\l>2$, as $H_1(
\Lambda_\l\otimes_{\Z F_{\l}} C_\bullet(F_{\l}))$ is not a
free module.

\example{Example 5.9} Consider the case $\l=3$, $n=5$, $m=1$.
Denote the generators of $\Lambda_4\otimes_{\Z P_{5,3}} C_2(P_{5,3}) =
(\Lambda_4)^{12}$ by $\{e_1,\dots,e_{12}\}$, and write $t_r = t_{r,5}$.
The elements
$$\matrix t_4(t_2-1)e_1+(1-t_1t_4)e_2+(t_2-1)e_4,\hfill&
(1-t_3)e_2+(t_2-1)e_3,\hfill\\
t_4(t_3-1)e_5+(1-t_1)e_7+(1-t_1)(t_3-1)e_8,&
t_4(t_3-1)e_6+(1-t_2t_4)e_7+(t_3-1)e_8,\hfill\\
(1-t_2)e_9+(t_1-1)e_{10},\hfill&
(1-t_3t_4)e_9+(t_1-1)e_{11}+t_3(t_1-1)e_{12},
\endmatrix
$$ generate a free, rank 6 submodule $M$ of
$\Lambda_4\otimes_{\Z P_{5,3}}C_2(P_{5,3})$.  Checking that this
submodule is invariant under the action of the representation
$\theta^1_{3,2,2}:P_3\to\GL(12,\Lambda_4)$, we obtain a
subrepresentation $\hat\theta^1_{3,2,2}:P_3\to\GL(6,\Lambda_4)$,
given by
$$\allowdisplaybreaks
\align
A_{1,2}&\mapsto
\pmatrix
t_1t_2t_4&0&0&0&0&0\\
t_3-1&1&0&0&0&0\\
0&0&1-t_1+t_1t_2t_4&t_1(1-t_1)&0&0\\
0&0&1-t_2t_4&t_1&0&0\\
0&0&0&0&t_1t_2&0\\
0&0&0&0&t_1(t_3t_4-1)&1
\endpmatrix,\\
\\
A_{1,3}&\mapsto
\pmatrix
t_3&1-t_1t_4&0&0&0&0\\
t_3(1-t_3)&1-t_3+t_1t_3t_4&0&0&0&0\\
0&0&t_1t_3&0&0&0\\
0&0&t_3(t_2t_4-1)&1&0&0\\
0&0&0&0&1&t_2-1\\
0&0&0&0&0&t_1t_3t_4
\endpmatrix,\\
\\
A_{2,3}&\mapsto
\pmatrix
1&t_2(t_1t_4-1)&0&0&0&0\\
0&t_2t_3&0&0&0&0\\
0&0&1&t_1-1&0&0\\
0&0&0&t_2t_3t_4&0&0\\
0&0&0&0&1-t_2+t_2t_3t_4&t_2(1-t_2)\\
0&0&0&0&1-t_3t_4&t_2
\endpmatrix.\\
\endalign$$
If one sets all $t_i$ equal to $t$ in the above representation
(i.e.~applies the map $\tilde\phi:\Lambda_4\to\Lambda$ defined by
$\phi(t_i)=t$), the resulting representation
$P_3\to\GL(6,\Lambda)$ is merely the restriction to $P_3$ of the
representation $\bar\beta_{3,2,2}^1:B_3\to\GL(6,\Lambda)$
discussed in Example~5.2.

The submodule $M\cong (\Lambda_4)^6$ considered above is {\it not}
isomorphic to the homology group $H_2(\Lambda_4 \otimes_{\Z P_{5,3}}
C_\bullet(P_{5,3}))$, which is not a free $\Lambda_4$-module.
This discrepancy disappears if one works over the complex numbers.
Let $\nu:\Z^4\to\C^*$ be a complex representation, and
$M\otimes_{\Z^4}\C_{\nu} \cong \C^6$ the corresponding vector space.
If the complex parameters $\nu(t_r)$ are sufficiently generic
(simply all different from $1$ in this instance), and
$\bar{\nu}(t_r):=\nu(t_r)^{-1}$, this vector space {\it is} isomorphic
to $H_2((\Lambda_4\otimes_{\Z P_{5,3}} C_\bullet(P_{5,3}))
\otimes_{\Z^4}\C_{\nu})\cong H_2(P_{5,3};\C_{\bar\nu}) \cong \C^6$
(see 4.12 for the first isomorphism). Similar genericity conditions are
addressed in detail in section~6.

Notice that the representation $\hat\theta_{3,2,2}^1$ decomposes into
a sum of 3 rank 2 representations of $P_3$.  More generally, one can
find a free, rank $\l(\l-1)$ submodule of
$\Lambda_{\l+1}\otimes_{\Z P_{\l+2,\l}}C_2(P_{\l+2,\l})$ which is
invariant under the action of the representation
$\theta_{\l,2,2}^1:P_{\l}\to\GL(\l(\l+1),\Lambda_{\l+1})$.  The
resulting subrepresentation $\hat\theta_{\l,2,2}^1$ decomposes into a
sum of $\l$ rank $\l-1$ representations of $P_{\l}$ for any $\l$, see
\cite{12}.

\endexample

\head 6. A Vanishing Theorem
\endhead

Let $P_\l$ denote the pure braid group on $\l$ strings.
In this section, we carry out a detailed analysis of the representations
which dictate the structure of the boundary maps of the free resolution
$C_\bullet(P_\l)$ of $\Z$ over the integral group ring $R=\Z P_\l$ given
by the construction of section 2.

\subhead{6.1}\endsubhead
The boundary maps of the chain complex
$C_\bullet(P_\l)$ are comprised of homomorphisms of the form
$\D_{p_1,\dots,p_k}:R^{p_1\cdots p_k}\to
R^{p_2\cdots p_k}$ defined
recursively by $\D_{p_1,\dots,p_k} =
-\tilde\rho_{p_k}(\D_{p_1\dots,p_{k-1}})$,
where $1\le p_1<p_2<\dots<p_k\le \l-1$.  Since
$\D_{p} = \left(
A_{1,p+1}-1\ \cdots\ A_{p,p+1}-1\right)^{\top}$, we have
$$\D_{p_1,\dots,p_k}= (-1)^{k-1}{\pmatrix
\phi(A_{1,p_1+1})-I&\cdots&\phi(A_{p_1,p_1+1})-I\endpmatrix}^{\top},$$
where $\phi = \tilde\rho_{p_k}\circ\cdots\circ\tilde\rho_{p_2}$ is the
composition of the representations
$\tilde\rho_{p_j}$, $2\le j\le k$, and $I$ is the identity matrix
of size $p_2\cdots p_k$.  Thus, to understand the chain complex
$C_\bullet(P_\l)$, we must make sense of the iterated compositions
of the representations $\tilde\rho_{p_j}$.

We first consider a single representation $\rho_p$.  To simplify
notation, we assume that
$p = \l-1$ and $r<s<\l$.  This representation is (essentially) given
by the Jacobian matrix of Fox derivatives of the action of $A_{r,s}$ on
the free group $F_{\l-1} = \langle
A_{1,\l},\dots,A_{\l-1,\l}\rangle$ exhibited in Example~1.3:
$$\rho_{\l-1}(A_{r,s}) = A_{r,s}\cdot\left({\partial A_{r,s}^{-1}
A_{i,\l} A_{r,s} \over \partial A_{j,\l}}\right).$$
For our immediate purposes (see below), rather that using the
``standard'' presentation of the pure braid group given in
Example~1.3, it is useful to work with a different generating
set for the free group $F_{\l -1}$.  Let
$$y_i=\cases
A_{i,\ell}&\text{if $i\le r$ or $i>s$,}\\
A_{r,\ell}A_{s,\ell}&\text{if $i=s$,}\\
A_{r,\ell}A_{i,\ell}A_{r,\ell}^{-1}&\text{if $r<i<s$.}\\
\endcases$$

A computation yields:
\proclaim{Lemma 6.2} The action of $A_{r,s}$ on the free
group $F_{\l-1}=\langle y_1,\dots,y_{\l-1}\rangle$ is given by
$$A_{r,s}^{-1} y_i A_{r,s} =
\cases y_i&\text{if $i\neq r$,}\\
y_s y_r y_s^{-1}&\text{if $i=r$.}\\
\endcases$$
\endproclaim
Computing Fox derivatives, we obtain:
\proclaim{Proposition 6.3} With respect to the generating set
$\{y_1,\dots,y_{\l-1}\}$ of $F_{\l-1}$, the matrix $U=(u_{i,j})$
of $\rho_{\l-1}(A_{r,s})$ is upper triangular with non-zero entries
$u_{r,r}=A_{r,s}A_{r,\l}A_{s,\l}=A_{r,s} y_s$,
$u_{r,s}=(1-y_r)A_{r,s}$, and $u_{k,k}=A_{r,s}$
for $k\ne r$.
\endproclaim

\subhead{6.4}\endsubhead
The above result may be viewed as a first step in an inductive
description of the matrix of
$\tilde\rho_{p_k}\circ\cdots\circ\tilde\rho_{p_1}(A_{r,s})$.
For example, with respect to the generating set
$\{A_{1,\l},\dots,A_{\l-1,\l}\}$ of $F_{\l-1}$, the matrix
of
$\tilde\rho_{\l-1}\circ\tilde\rho_{k-1}(A_{r,s})$ is, by Proposition~6.3,
similar to the matrix obtained by applying $\rho_{\l-1}$
to the entries of the matrix
$$\pmatrix A_{r,s}\cdot I_{r-1}&0&0\\ 0&U_{r,s}&0\\0&0&A_{r,s}\cdot
I_{k-s-1}\\
\endpmatrix,$$ where $U_{r,s}$ is an $(s-r+1)\times(s-r+1)$ upper
triangular matrix with diagonal entries
$A_{r,s}A_{r,k}A_{s,k},A_{r,s}, \dots, A_{r,s}$.
Thus, to understand the composition
$\tilde\rho_{\l-1}\circ\tilde\rho_{k-1}(A_{r,s})$,
we must understand not only
$\rho_{\l-1}(A_{r,s})$, but $\rho_{\l-1}(A_{r,s}A_{r,k}A_{s,k})$
as well.

Let $J = (j_1,j_2,\dots,j_p)$, where $1\le j_1<j_2<\cdots<j_p\le
\l$, and consider the pure braid $A_J$ defined by
$$A_J = (A_{j_1,j_2})(A_{j_1,j_3}A_{j_2,j_3})(A_{j_1,j_4}A_{j_2,j_4}
A_{j_3,j_4})
\cdots\cdots(A_{j_1,j_p}\cdots A_{j_{p-1},j_p}).$$
If $j_p \le \l-1$, then $A_J$ acts on the free group $F_{\l-1} =
\langle A_{1,\l},\dots,A_{\l-1,\l}\rangle$ by conjugation.
We once again work with a different generating set of the free
group $F_{\l-1}$.  Let
$$y_i=\cases
A_{i,\ell}&\text{if $i\le j_1$ or $i>j_p$,}\\
A_{j_1,\ell}\cdots A_{j_k,\ell}&\text{if $i= j_k$, $2\le k\le p$,}\\
(A_{j_1,\ell}\cdots
A_{j_k,\ell})A_{i,\ell}(A_{j_1,\ell}\cdots A_{j_k,\ell})^{-1}
&\text{if $j_k<i<j_{k+1}$, $1\le k\le p-1$.}\\
\endcases$$

\proclaim{Lemma 6.5} The action of $A_J$ on the free group
$F_{\l-1}=\langle y_1,\dots,y_{\l-1}\rangle$ is given by
$$A_J^{-1} y_i A_J =
\cases y_i&\text{if $i\neq j_k$, $1\le k\le p-1$,}\\
y_{j_p} y_i y_{j_p}^{-1}&\text{if $i = j_k$, $1\le k\le p-1$.}\\
\endcases
$$
\endproclaim
\demo{Proof} If $p=2$, the above is merely the reformulation of
the defining relations of the pure braid group $P_\l$
considered in Lemma~6.2.  Using induction on $p$ and the
defining relations in the pure braid group, one checks that the
action of $A_J$ on $F_{\l-1}$ is as asserted. \quad\qed
\enddemo
Computing Fox derivatives once again, we obtain:
\proclaim{Proposition 6.6} With respect to the generating set
$\{y_1,\dots,y_{\l-1}\}$ of $F_{\l-1}$, the
matrix $U= (u_{i,j})$ of
$\rho_{\l-1}(A_J)$ is upper triangular with non-zero entries
$u_{j_k,j_k} = A_J y_{j_p}$ and $u_{j_k,j_p}=(1-y_{j_k})A_J$
for $1\le k\le p-1$, and $u_{i,i} = A_J$ for $i \neq j_k$, $1\le
k\le p-1$.
\endproclaim

\subhead{6.7}\endsubhead We now use these results to study the
homology of the pure braid group (and related groups) with
non-trivial coefficients.  Recall that, for $1\le \l < n$,
$P_{n,\l}=\ker (P_n\to P_\l)$  denotes the kernel of the
homomorphism given by $A_{r,s}\mapsto 1$ if
$s\le\l$ and  $A_{r,s}\mapsto A_{r,s}$ if $s>\l$.
We consider complex representations,
$\nu:P_{n,\l}\to \GL(m,\C)$, of these groups.
For $J=(j_1,\dots,j_p)$ such that $A_J = (A_{j_1,j_2})\cdot
(A_{j_1,j_3}A_{j_2,j_3})\cdots (A_{j_1,j_p}
\cdots A_{j_{p-1},j_p})$ is an element of $P_{n,\l}$, write
$\nu(A_J)=v_{_J}$. Let $V$ denote the $(\Z P_{n,\l})$-module
(resp.~local coefficient system on $K(P_{n,\l},1)$) corresponding to
$\nu$.

\definition{Definition 6.8} Fix $q\ge 0$.  The
representation $\nu:P_{n,\l}\to\GL(m,\C)$ is said to be {\it
quasi-generic through rank $q$} if, for each $J$ such that $2\le
|J|\le q+2$ and $A_J\in P_{n,\l}$, the eigenvalues of the
matrix $v_{_J}$ are all different from $1$.
\enddefinition

\remark{Remark 6.9} If $\nu$ is quasi-generic through rank $q$,
repeated application of Proposition 6.6 shows that the
eigenvalues of the matrix of $\tilde\nu\circ\phi(A_{i,j})$ are all
different from $1$, where $\tilde\nu:\Z P_{n,\l}\to \End(\C^m)$
denotes the linear extension of $\nu$,
$\phi=\tilde\rho_{p_m}\circ\cdots\circ\tilde\rho_{p_1}$, and
$m\le q+1$.  Thus, under these conditions, the endomorphism
$I-\tilde\nu\circ\phi(A_{i,j})$ is in fact an isomorphism.
This observation motivates the following, which generalizes the
vanishing theorem found in \cite{29} (see also~\cite{46}, \cite{50}).
\endremark

\proclaim{Theorem 6.10} If $\nu:P_{n,\l}\to\GL(m,\C)$ is
quasi-generic through rank $q$, then the homology
groups of $P_{n,\l}$ with coefficients in $V$
vanish for $0\le i\le \min\{q,n-\l-1\}$.
\endproclaim
\demo{Proof}  The proof is by induction on the cohomological dimension,
$d=\text{cd}(P_{n,\l})=n-\l$, of the group $P_{n,\l}$.

If $d=n-\l=1$, we have $P_{n,\l}=F_{n-1}$ and
$\min\{q,n-\l-1\}=0$.  In this instance, the hypothesis of the
theorem merely states that the eigenvalues of the matrices
$\nu(A_{i,\l})$ are different from $1$, and it follows easily
that $H_0(F_{n-1};V)=0$.

In the general case, write $C_\bullet(n,\l)=C_\bullet(P_{n,\l})
\otimes_{ P_{n,\l}} V$, and denote the boundary maps of this complex
by $\partial_i(n,\l)$.  The restriction of the representation
$\nu$ to the subgroup $P_{n,\l+1}$ of $P_{n,\l}$ gives rise to a
$(\Z P_{n,\l+1})$-module which we continue to denote by $V$.
By induction, we have $H_i(P_{n,\l+1};V)=0$ for $0\le i\le
\min\{q,n-\l-2\}$.

As in the proof of Theorem~2.11, we exploit
the fact that the complex $C_\bullet(P_{n,\l})$ may be realized as
the mapping cone of $\Xi_\bullet:[C_\bullet(n,\l+1)]^\l\to
C_\bullet(n,\l+1)$, where
$C_\bullet(n,\l+1)=C_\bullet(P_{n,\l+1})\otimes_{P_{n,\l+1}}V$.
Thus we have a short exact sequence of chain complexes
$$0 \to C_\bullet(n,\l+1) \to C_\bullet(n,\l) \to
[C_{\bullet-1}(n,\l+1)]^\l \to 0.$$
It follows immediately from the corresponding long exact sequence
in homology:
$$\dots\to [H_i(P_{n,\l+1};V)]^\l \to H_i(P_{n,\l+1};V) \to
H_i(P_{n,\l};V) \to [H_{i-1}(P_{n,\l+1};V)]^\l \to \dots$$
that $H_i(P_{n,\l};V)=0$ for $0\le i\le
\min\{q,n-\l-2\}$.  This completes the proof if $q\le n-\l-2$.

If $q\ge n-\l-1$,  since we have $H_{i}(P_{n,\l};V)=0$ for
$i\le n-\l-2$ by the previous paragraph, it remains to show that
$H_{n-\l-1}(P_{n,\l};V)=0$.  For this, it suffices to show that
$\dim\im\partial_{n-\l}(n,\l)=\dim\ker\partial_{n-\l-1}(n,\l)$.
The boundary map $\partial_{n-\l}(n,\l)$ is of the form
$$\partial_{n-\l}(n,\l) = \pmatrix \D_{\l,\l+1,\dots,n-1} &
[\partial_{n-\l-1}(n,\l+1)]^\l\endpmatrix.$$
The map $\D_{\l,\l+1,\dots,n-1}$ is given by
$$\D_{\l,\l+1,\dots,n-1}=(-1)^{n-\l}{\pmatrix
\tilde\nu\circ\phi(A_{1,\l+1})-I&\cdots&
\tilde\nu\circ\phi(A_{\l,\l+1})-I\endpmatrix}^{\top},$$
where $\phi = \tilde\rho_{n-1}\circ\cdots\circ\tilde\rho_{\l+1}$
is the composition of the representations $\tilde\rho_{j}$,
$\l+1\le j\le n-1$, and $I$ is the
identity matrix of size $m\cdot(\l+1)\cdots (n-1)$.  Thus
$\partial_{n-\l}(n,\l)$ has a submatrix of the form
$$\pmatrix
(-1)^{n-\l}(\tilde\nu\circ\phi(A_{1,\l+1})-I) &
0\\ *& [\partial_{n-\l-1}(n,\l+1)]^{\l-1}
\endpmatrix,$$
and
$$\dim\im(\partial_{n-\l}(n,\l))\ge (\l-1)\cdot
\dim\im(\partial_{n-\l-1}(n,\l+1)) +
\rank(\tilde\nu\circ\phi(A_{1,\l+1})-I).$$

Now the conditions on the representation $\nu$ assure that
$\tilde\nu\circ\phi(A_{1,\l+1})-I$ is an invertible matrix,
hence has rank
$m\cdot(\l+1)\cdots (n-1)$.  Since $H_i(P_{n,\l};V)=0$ for
$i\le n-\l-2$, we compute $\dim\ker\partial_{n-\l-1}(n,\l)$
using an Euler characteristic argument.  It then follows easily that
$$(\l-1)\cdot\dim\im(\partial_{n-\l-1}(n,\l+1)) +
\rank(\tilde\nu\circ\phi(A_{1,\l+1})-I) =
\dim\ker(\partial_{n-\l-1}(n,\l)).\quad\qed$$
\enddemo

\proclaim{Corollary 6.11} If $q \ge n-\l-1$ and
$\nu:P_{n,\l}\to\GL(m,\C)$ is quasi-generic through
rank $q$, then the homology group
$H_{n-\l}(P_{n,\l};V)$ has rank $m\cdot(n-2)!/(\l-2)!$ if
$\l\ge 2$, and is trivial if $\l =1$.
\endproclaim

\remark{Remark 6.12} Note that if $\nu:P_{n,1}\to\GL(m,\C)$ is a
quasi-generic representation through rank $q$ of the (entire)
pure braid group, and $q\ge n-2$, then all homology groups
$H_{i}(P_{n,1};V)$ vanish.  That is, the chain complex
$C_\bullet(n,1) = C_\bullet(P_n)\otimes_{ P_n} V$ is
acyclic.
\endremark

\remark{Remark 6.13} For certain rank one local coefficient
systems $V$ on $K(P_{n,\l},1)$, Schechtman and Varchenko \cite{47}
use cycles in $H_{n-\l}(K(P_{n,\l},1);V)$, together with
differential forms in the (de~Rham) cohomology group
$H^{n-\l}(K(P_{n,\l},1);{V}^*)$
with coefficients in the dual local system, to generate
solutions of the Knizhnik-Zamolodchikov equations in terms
of generalized hypergeometric functions (see also \cite{16},
\cite{31}).  In the instances where the local coefficient systems
constructed in \cite{47} arise from representations which are
quasi-generic through rank $n-\l-1$, it follows from
Corollary~6.11 that there are $(n-2)!/(\l-2)!$ linearly
independent solutions of the corresponding KZ equations.
\endremark

\head 7. Milnor fibrations
\endhead

In this section, we discuss how the chain complex constructed in
section~2 and the vanishing theorem of the previous section may be used
in the study of Milnor fibrations.

\subhead{7.1}\endsubhead   Let $\A$ be a central fiber-type
arrangement in $\C^\l$ with complement $M=M(\A)$ and
exponents $\{1=d_1,d_2,\dots,d_\l\}$ (see \cite{18},
\cite{43}).  Then the fundamental group $G$ of
$M$ may be realized as an iterated semidirect product of free
groups, $G\cong F_{d_\l}\rtimes\cdots\rtimes F_{d_1}$, and $M$
is an Eilenberg-MacLane space of type $K(G,1)$. Let
$Q=Q(\A)$ be a defining polynomial of $\A$.  Then, since $\A$ is central,
$Q$ is homogeneous of degree $n=\sum d_q=|\A|$ and we have a (global)
Milnor fibration $Q: M\to\C^*$, with fiber $F=Q^{-1}(1)$, and
monodromy $h:F\to F$ given by multiplication by $\xi=\exp(2\pi\i/n)$,
\cite{38}.  The Milnor fiber $F=F(\A)$ has the homotopy type of
an $(\l-1)$-dimensional $K(\pi,1)$ space, where $\pi=\pi_1(F)$.

Since $F$, $M$, and $\C^*$ are Eilenberg-MacLane spaces, the
homotopy exact sequence of the Milnor fibration reduces to
$1\to\pi\to G\rightleftarrows\Z\to 1$.  Therefore, by
Shapiro's Lemma (see \cite{9}), we have
$$H_*(F;\Z)=H_*(\pi;\Z)=H_*(G;\Z G\otimes_{\Z\pi}\Z)=H_*(G;\Z\Z).$$
Thus the construction of section~2 provides an algorithm for
computing the integral homology of the Milnor fiber of an arbitrary
fiber-type arrangement.  This algorithm also applies to certain
non-linearly fibered arrangements, such as the Coxeter arrangements
of type $\operatorname{D}_\l$ mentioned in Example 1.8.

\subhead{7.2}\endsubhead
We have carried out this computation for the braid arrangements
$\A_\l$, for $\l\le 5$.  The results are tabulated below.
For $\l\le 4$, these explicit results have been obtained by other
means (see e.g.~\cite{11} and the references therein).
To the best of our knowledge, the results for $\A_5$ were
previously unknown.  Since there is no torsion in the homology
of these Milnor fibers, we list only the Betti numbers.
We conjecture that the homology of the Milnor fiber of the braid
arrangement $\A_\l$ is torsion free for any $\l$.
$$
\centerline{
\vbox{\offinterlineskip
\hrule
\halign{&\vrule#&\strut\quad\hfil#\hfil\quad\cr
height2pt&\omit&&\omit&\cr&\hfil&&$b_0(F)$
\hfil&&\hfil$b_1(F)$\hfil&&\hfil$b_2(F)$\hfil&&\hfil$b_3(F)$\hfil&\cr
height2pt&\omit&&\omit&&\omit&&\omit&&\omit&\cr
\noalign{\hrule}
height2pt&\omit&&\omit&&\omit&&\omit&&\omit&\cr
&$\A_3$&&1&&4&&&&&\cr
height2pt&\omit&&\omit&&\omit&&\omit&&\omit&\cr
\noalign{\hrule}
height2pt&\omit&&\omit&&\omit&&\omit&&\omit&\cr
&$\A_4$&&1&&7&&18&&&\cr
height2pt&\omit&&\omit&&\omit&&\omit&&\omit&\cr
\noalign{\hrule}
height2pt&\omit&&\omit&&\omit&&\omit&&\omit&\cr
&$\A_5$&&1&&9&&28&&80&\cr
height2pt&\omit&&\omit&&\omit&&\omit&&\omit&\cr
\noalign{\hrule}}}}
$$

For arbitrary $\l$, we can compute the Euler characteristic
of $F(\A_\l)$.  Indeed, the Milnor fiber of $\A_\l$ is a cyclic
$\binom{\l}{2}$-sheeted cover of the complement of the
projectivized braid arrangement, $M(\A_{\l}^*)$, (see e.g.~\cite{43}).
Since $M(\A_{\l}^*)$ is a $K(P_{\l,2},1)$ space, it follows from the
LCS~formula that $\chi(M(\A_{\l}^*))=(-1)^{\l}(\l-2)!$.  Thus
$\DS{\chi(F(\A_\l))=(-1)^{\l}{\l}!/2}$.

\subhead{7.3}\endsubhead
The construction of section~2 may also be used to compute the homology
eigenspaces of the algebraic monodromy of the Milnor fibration.
In \cite{11}, it is shown that, for an arbitrary central arrangement
$\A$, the $\xi^k$-eigenspace of the monodromy is isomorphic to
$H_*(M(\A^*);{V}_k)$, the homology of the complement of the
projectivization of $\A$ with coefficients in a complex rank one local
system ${V}_k$.  This local system is induced by the representation
$\nu_k:\pi_1(M(\A^*))\to\C^*$ which sends each meridian of $\A^*$ to
$\xi^k$.

We have computed the homology eigenspaces of the monodromy of the Milnor
fibration of the braid arrangement
$\A_\l$, for $\l\le 5$.  The characteristic polynomials,
$p_i(t)$, of the maps induced in homology by the monodromy,
$h_*:H_i(F)\to H_i(F)$, are given below.
$$
\centerline{
\vbox{\offinterlineskip
\hrule
\halign{&\vrule#&\strut\ \,\hfil#\hfil\ \,\cr
height2pt&\omit&&\omit&\cr&\hfil&&$p_0(t)$
\hfil&&\hfil$p_1(t)$\hfil&&\hfil$p_2(t)$\hfil&&\hfil$p_3(t)$\hfil&\cr
height2pt&\omit&&\omit&&\omit&&\omit&&\omit&\cr
\noalign{\hrule}
height2pt&\omit&&\omit&&\omit&&\omit&&\omit&\cr
&$\A_3$&&$1-t$&&$(1-t)(1-t^3)$&&&&&\cr
height2pt&\omit&&\omit&&\omit&&\omit&&\omit&\cr
\noalign{\hrule}
height2pt&\omit&&\omit&&\omit&&\omit&&\omit&\cr
&$\A_4$&&$1-t$&&$(1-t)^4(1-t^3)$&&$(1-t)^3(1-t^3)(1-t^6)^2$&&&\cr
height2pt&\omit&&\omit&&\omit&&\omit&&\omit&\cr
\noalign{\hrule}
height2pt&\omit&&\omit&&\omit&&\omit&&\omit&\cr
&$\A_5$&&$1-t$&&$(1-t)^9$&&$(1-t)^{26}(1-t^2)$&&
$(1-t)^{18}(1-t^2)(1-t^{10})^6$&\cr
height2pt&\omit&&\omit&&\omit&&\omit&&\omit&\cr
\noalign{\hrule}}}}
$$
For arbitrary $\l$, the zeta function of the monodromy
is given by
$$\zeta(t)=p_0(t)^{-1}\cdot p_1(t)\cdot p_2(t)^{-1}\cdot p_3(t)\cdots
p_{\l-2}(t)^{\pm 1} = \left(1-t^n\right)^{(-1)^{\l+1}(\l-2)!},$$
where $n=\binom{\l}{2}$, see \cite{38}.

\remark{Remark 7.4} Since the complement of the projectivized
braid arrangement $M(\A_{\l}^*)$ is a $K(P_{\l,2},1)$ space, Theorem~6.10
may be used to obtain partial results on the homology eigenspaces of the
monodromy of the Milnor fibration of the braid arrangement for general
$\l$. For instance, if $n=\binom{\l}{2}$ is not divisible by $2$ or $3$,
then for
$1\le k \le n$, the rank one representation $\nu_k$ of
$\pi_1(M(\A_{\l}^*))=P_{\l,2}$ is quasi-generic through rank $2$. It
follows from Theorem~6.10 (see \cite{11} for details) that
$H_i(F(\A_{\l}))=H_i(M(\A_{\l}^*))$ for $i\le 2$.

For arbitrary $\l$, if $\xi^k$ is a primitive $n^{\text{th}}$ root of
unity, then the representation $\nu_k$
is quasi-generic through rank $\l-3$.
Thus the $\xi^k$-eigenspace of the monodromy is ``concentrated'' in
dimension $\l-2$, that is $H_i(M(\A_{\l}^*);{V}_k)=0$ if $i \neq
\l-2$.
\endremark

\subhead{7.5}\endsubhead  The approach outlined above
also gives information on the Milnor fibration
of the discriminant
singularity $\dd_\l$ in $\C^\l$.  (See \cite{44} for another
possible way to attack this problem.)
As noticed by Arnol'd, the complement, $M(\dd_\l)$,
is the configuration space of the set of $\l$ (unordered) points
in $\C$, and thus is a $K(B_\l,1)$-space, see e.g. \cite{23}.
The Milnor fibration $M(\dd_\l)\to\C^*$ induces on $\pi_1$
the abelianization map $\ab: B_\l\to\Z$, see \cite{22}.
Hence, the Milnor fiber, $F(\dd_\l)$, is a $K(B'_\l,1)$-space.

For arbitrary $\l$ we can compute the Euler characteristic of
$F(\dd_\l)$.  Indeed, the usual
symmetric group covering $M(\A_\l) \to M(\dd_\l)$ restricts on Milnor
fibers to an
alternating group covering $F(\A_\l) \to F(\dd_\l)$.  Thus
$$\chi(F(\dd_\l))=\frac{(-1)^{\l}{\l}!/2}{{\l}!/2} =(-1)^{\l}.$$
However, computing the homology of $F(\dd_\l)$ is, in general, a
substantially harder task.

The case $\l=3$ is well-known:  $\dd_3$ is the product of the
cusp singularity with $\C$.  Hence, $\pi_1(F(\dd_3))=F_2$ and the
Milnor number, $b_1(F(\dd_3))$, is 2.

The case $\l=4$ is not as well-known.  In \cite{36}, Massey uses
L\^ e numbers to find upper bounds for the Betti numbers of the Milnor
fiber of $\dd_4$, the product of the swallowtail singularity with $\C$.
He finds
$b_1(F(\dd_4))=b_2(F(\dd_4))\le 5$.
Our approach gives the exact answer.
Let $G:=B_4'=\pi_1(F(\dd_4))$.  It is shown in \cite{22} that $G$ is a
semidirect product $G_2\rtimes G_1$.  The action of
$G_1=F_2=\langle x_{1,1}, x_{2,1}\rangle$ on
$G_2=F_2=\langle x_{1,2}, x_{2,2}\rangle$ is given by
$$
x_{1,1}:
\cases
x_{1,2}\mapsto x_{1,2}x_{2,2}^{-1} x_{1,2}^2\\
x_{2,2}\mapsto x_{1,2}
\endcases
\qquad
x_{2,1}:
\cases
x_{1,2}\mapsto x_{1,2}x_{2,2}^{-1} x_{1,2}^3\\
x_{2,2}\mapsto x_{1,2}x_{2,2}^{-1} x_{1,2}^4
\endcases
$$
A computation reveals that $H_1(G)=\Z^2, H_2(G)=\Z^2$.  Thus
$b_1(F(\dd_4))=b_2(F(\dd_4))=2$.

\enddocument